\documentclass{au}
\usepackage{amssymb}

\newcommand{\<}{\kern.0833em}

\newtheorem{theorem}{Theorem}
\newtheorem{lemma}[theorem]{Lemma}
\newtheorem{corollary}[theorem]{Corollary}
\newtheorem{proposition}[theorem]{Proposition}
\newtheorem{definition}[theorem]{Definition}
\newtheorem{question}[theorem]{Question}

\newcommand{\Z}{\mathbb Z}

\newcommand{\Sym}{\mathrm{Sym}}
\newcommand{\strt}{\rule{0em}{.6em}} 
\newcommand{\cl}{\mathrm{cl\<}}
\newcommand{\preq}{\preccurlyeq}
\newcommand{\sucq}{\succcurlyeq}
\newcommand{\eps}{\varepsilon}
\newcommand{\cj}{^{\rm cj}}
\renewcommand{\leq}{\leqslant}
\renewcommand{\geq}{\geqslant}

\makeatletter
\newcommand{\xlabel}{\stepcounter{equation}
  \gdef\@currentlabel{\p@equation\theequation}{\rm(\@currentlabel)}}
\makeatother
\newenvironment{xlist}
  {\begin{list}{\xlabel}
    {\setlength{\rightmargin}{20pt}
     \setlength{\leftmargin}{37pt}
     \setlength{\labelsep}{20pt}
     \setlength{\labelwidth}{20pt}}}
  {\end{list}}

\begin{document}

\title{Closed subgroups of the infinite symmetric group}
\thanks{
Second author's research supported by the United States--Israel
Binational Science Foundation.
\protect\\
\indent
A preprint version of this paper is readable at
http://math.berkeley.edu/%
{$\!\sim\!$}gbergman/\linebreak[3]papers/%
Sym\<\_\,Omega:2.\{tex,dvi\}\<,
at http://shelah.logic.at
as publication 823, and at
arXiv:math\linebreak[3].GR/0401305\,.
}

\author{George M. Bergman}
\address[G. Bergman]{University of California\\
Berkeley, CA 94720-3840, USA}
\email{gbergman@math.berkeley.edu}

\author{Saharon Shelah}
\address[S. Shelah]{University of Jerusalem\\
Einstein Institute of Mathematics\\
Jerusalem 91904, Israel\\
and
Hill Center-Busch Campus\\
Rutgers, The State University of New Jersey\\
Piscataway, NJ 08854-8019, USA}
\email{shelah@math.huji.ac.il}

\begin{abstract}
Let $S=\Sym(\Omega)$ be the group of all permutations of a countably
infinite set $\Omega,$ and for subgroups $G_1,\,G_2\leq S$ let us
write $G_1\approx G_2$ if there exists a finite set $U\subseteq S$
such that $\langle\<G_1\cup U\<\rangle = \langle\<G_2\cup U\<\rangle.$
It is shown that the subgroups closed in the function topology on $S$
lie in precisely four equivalence classes under this relation.
Which of these classes a closed subgroup $G$ belongs to depends on
which of the following statements about pointwise stabilizer subgroups
$G_{(\Gamma)}$ of finite subsets $\Gamma\subseteq\Omega$ holds:
\\[2pt]
(i)~For every finite set $\Gamma,$ the subgroup $G_{(\Gamma)}$ has
at least one infinite orbit in $\Omega\<.$
\\[2pt]
(ii)~There exist finite
sets $\Gamma$ such that all orbits of $G_{(\Gamma)}$
are finite, but none such that the cardinalities
of these orbits have a common finite bound.
\\[2pt]
(iii)~There exist finite sets $\Gamma$ such that the cardinalities
of the orbits of $G_{(\Gamma)}$ have a common finite bound,
but none such that $G_{(\Gamma)}=\{1\}.$
\\[2pt]
(iv)~There exist finite sets $\Gamma$ such that $G_{(\Gamma)}=\{1\}.$
\vspace{2pt}

Some related results and topics for further investigation are noted.
\end{abstract}
\dedicatory{In honor of Walter Taylor, on his not-yet-retirement}
\maketitle

\section{Introduction.}\label{S.Intro}
In \cite[Theorem~1.1]{DM&PN}, Macpherson and Neumann show
that for $\Omega$ an infinite set, the group $S=\Sym(\Omega)$ is not
the union of a chain of $\leq|\<\Omega\<|$ proper subgroups.
It follows that if $G$ is a subgroup of $S,$ and if
$S=\langle\<G\cup U\rangle$
for some set $U\subseteq S$ of cardinality $\leq|\<\Omega\<|,$
then one may replace $U$ by a finite subset of $U$ in this equation.
Galvin~\cite{FG} has shown that in this situation one can even replace
$U$ by a singleton, though not necessarily one contained in $U$ or
even in $\langle\<U\rangle.$

Thus we have a wide gap -- between subgroups $G$ over which it
is ``easy'' to generate $S$ (where one additional element will
do), and all others, over which it is ``hard'' (even $|\<\Omega\<|$
elements will not suffice).
It is natural to wonder how one can tell to which
sort a given subgroup belongs.
There is probably no simple answer for arbitrary subgroups;
but we will show, for $\Omega$ countable, that if our subgroup is
closed in the function topology on $S,$ then one element suffices
if and only if $G$ satisfies condition~(i) of the above abstract.
The method of proof generalizes to give the four-way
classification of closed subgroups asserted there.

(The four conditions of that classification could be
stated more succinctly, if not so transparently to the
non-set-theorist, by writing $\lambda$ for the least cardinal
such that for some finite subset $\Gamma\subseteq\Omega,$ all orbits
of $G_{(\Gamma)}$ in $\Omega$ have cardinality $<\lambda\<.$
Then the conditions are (i)~$\lambda=\aleph_1,$
(ii)~$\lambda=\aleph_0,$ (iii)~$3\leq\lambda<\aleph_0$ and
(iv)~$\lambda=2.$
But we shall express them below in the more
mundane style of the abstract.)

The proofs of the above results will
occupy~\S\S\ref{S.Defs}--\ref{S.countable} of this note.
In~\S\S\ref{S.q.FN}--\ref{S.q.etc} we note some related observations,
questions, and possible directions for further investigation.

The authors are indebted to Peter Biryukov and Zachary Mesyan
for corrections to earlier drafts of this note.

\section{Definitions, conventions, and basic
observations.}\label{S.Defs}

As usual, ``$\leq$'' appearing before the symbol for a group means
``is a subgroup
of'', and $\langle\,...\,\rangle$ denotes ``the subgroup generated by''.

We take our notation on permutation groups from~\cite{DM&PN}.
Thus, if $\Omega$ is a set, $\Sym(\Omega)$ will
denote the group of all permutations of $\Omega,$
and such permutations will be written to the right of their arguments.
Given a subgroup $G\leq \Sym(\Omega)$ and a subset
$\Sigma\subseteq\Omega,$ the symbol $G_{(\Sigma)}$ will denote
the subgroup of elements of $G$ that stabilize $\Sigma$ pointwise, and
$G_{\{\Sigma\}}$ the larger subgroup $\{\<f\in G:\Sigma f=\Sigma\<\}.$
Elements of $\Omega$ will generally be denoted $\alpha,\,\beta,\ldots$.

The cardinality of a set $U$ will be denoted $|\<U|.$
Each cardinal is understood to be the least ordinal of its
cardinality, and each ordinal to be the set of all smaller ordinals.
The successor cardinal of a cardinal $\kappa$ is denoted~$\kappa^+.$
\vspace{6pt}

Let us now define, in greater generality than we did in the abstract,
the relations we will be studying.

\begin{definition}\label{D.approx}
If $S$ is a group, $\kappa$ an infinite cardinal, and $G_1,\,G_2$
subgroups of $S,$ we shall write $G_1\preq_{\strt\kappa,S} G_2$
if there exists a subset $U\subseteq S$ of cardinality $<\kappa$
such that $G_1\leq \langle\<G_2\<\cup\<U\rangle.$
If $G_1\preq_{\strt\kappa,S} G_2$ and
$G_2\preq_{\strt\kappa,S} G_1,$
we shall write $G_1\approx_{\strt\kappa,S} G_2,$
while if $G_1\preq_{\strt\kappa,S} G_2$ and
$G_2\not\preq_{\strt\kappa,S} G_1,$
we shall write $G_1\prec_{\strt\kappa,S} G_2.$

We will generally omit the subscript $S,$ and often
$\kappa$ as well, when their values are clear from the context.
\end{definition}

Clearly $\preq_{\strt\kappa,S}$ is a preorder on subgroups
of $S,$ hence $\approx_{\strt\kappa,S}$ is an
equivalence relation, equivalent to the
assertion that there exists $U\subseteq S$ of cardinality $<\kappa$
such that $\langle\<G_1\cup U\rangle=\langle\<G_2\cup U\rangle.$
Note that conjugate subgroups of $S$ are
$\approx_{\strt\kappa}\!$-equivalent for all $\kappa\<.$
If $G_1$ and $G_2$ are \mbox{$\approx_{\strt\kappa}\!$-equivalent,}
we see that they are \mbox{$\approx_{\strt\kappa}\!$-equivalent}
to $\langle\<G_1\cup G_2\<\rangle.$
(However they need not be
$\approx_{\strt\kappa}\!$-equivalent to $G_1\cap\<G_2.$
For instance, if $S=\Sym(\Z)$ and $G_1,\,G_2$
are the pointwise stabilizers of the sets of positive, respectively
negative integers, then they are conjugate, so
$G_1\approx_{\strt\aleph_0}G_2,$ but $G_1\cap G_2=\{1\},$ which is
not $\approx_{\strt\aleph_0}\!$-equivalent to $G_1$ and $G_2,$ since
the latter groups are uncountable, hence not finitely generated.)

We note

\begin{lemma}\label{L.G_*G}
Let $\Omega$ be an infinite set, $G$ a subgroup of $S=\Sym(\Omega),$
and $\Gamma\subseteq\Omega$ a subset such that
$|\<\Omega\<|^{|\Gamma|}\leq|\<\Omega\<|$
{\rm(}\!\<e.g., a finite subset{\rm)}.
Then $G_{(\Gamma)}\approx_{\strt|\<\Omega\<|^+} G.$

\end{lemma}\begin{proof}
Elements of $G$ in distinct right cosets of
$G_{(\Gamma)}$ have distinct behaviors on $\Gamma,$
hence if $R$ is a set of representatives of these
right cosets, $|R\<|\leq|\<\Omega\<|^{|\Gamma|}\leq|\<\Omega\<|.$
Clearly, $\langle\<G_{(\Gamma)}\cup R\<\rangle=G=
\langle\<G\cup R\<\rangle,$ so
$G_{(\Gamma)}\approx_{\strt|\<\Omega\<|^+} G,$ as claimed.
\end{proof}

Two results from the literature have important
consequences for these relations:
\begin{lemma}[{\rm\cite{FG}, \cite{DM&PN}}\textbf{}]\label{L._0<=>_1}
Let $\Omega$ be an infinite set.
Then on subgroups of $S=\Sym(\Omega),$\\[2pt]
{\rm(i)} The binary relation
$\preq_{\strt\aleph_0}$ coincides with $\preq_{\strt\aleph_1}$
{\rm(}hence $\approx_{\strt\aleph_0}$ coincides with
$\approx_{\strt\aleph_1}).$\\[2pt]
{\rm(ii)} The unary relation $\approx_{\strt\aleph_0}S$
coincides with $\approx_{\strt|\<\Omega\<|^{^+}}S.$
\end{lemma}\begin{proof}
(i) follows from~\cite[Theorem~3.3]{FG}, which says that every countably
generated subgroup of $S$ is contained in a $2\!\<$-generator subgroup.

We claim that~(ii) is a consequence of \cite[Theorem~1.1]{DM&PN}
$(=$\cite[Theorem~5]{Sym_Omega:1}) which, as recalled
in \S\ref{S.Intro}, says that any chain of proper subgroups of
$S$ having $S$ as union must have $>|\<\Omega\<|$ terms.
For if $G\approx_{\strt|\<\Omega\<|^{^+}}S,$ then among
subsets $U\subseteq S$ of cardinality $\leq|\<\Omega\<|$ such that
$\langle\<G\cup U\rangle=S,$ we can choose one of least cardinality.
Index this set $U$ as $\{g_i:i\in|\<U|\}.$
Then the subgroups $G_i=\langle\<G\cup \{g_j:j<i\}\<\rangle$
$(i\in|\<U|)$ form a chain of $\leq |\<\Omega\<|$ proper subgroups
of $S,$ which if $|\<U|$ were infinite would have
union $S,$ contradicting the result of \cite{DM&PN} quoted.
So $U$ is finite, so $G\approx_{\strt\aleph_0}S.$
\end{proof}

We have not introduced the symbols $\preq_{\strt\kappa}$ and
$\approx_{\strt\kappa}$ for $\kappa$ finite because in
general these relations are not transitive; rather, one has
$G_1\preq_{\<m} G_2\preq_{\<n} G_3\implies G_1\preq_{\<m+n-1}G_3.$
However, the proof of~(i) above shows that in
groups of the form $\Sym(\Omega),$
$\preq_{\strt\aleph_0}$ is equivalent to $\preq_3.$
Thus, in such groups, we do have transitivity of
$\preq_{\strt\kappa}$ and $\approx_{\strt\kappa}$ when
$3\leq \kappa<\aleph_0,$ but have no need for these symbols.
These observations are also the reason why in
\S\S\ref{S.inf_orbs}--\ref{S.bdd_orbs} we won't make
``stronger'' assertions than $\approx_{\aleph_0},$ though some of our
constructions will lead to sets $U$ of explicit finite cardinalities.

(Incidentally, \cite[Theorem~5.7]{FG} shows that on $\Sym(\Omega),$
the unary relation
$\approx_{\strt\aleph_0}S$ is even equivalent to $\approx_2 S.)$

\section{Generalized metrics.}\label{S.metrics}

In this section we will prove a result which, for
$\Omega$ countable, will imply that closed
subgroups of $\Sym(\Omega)$ falling under different cases
of the classification described in the abstract
are indeed $\approx_{\strt\aleph_0}\!$-inequivalent.

To motivate our approach, let us sketch
a quick proof that if $\Omega$ is any set such that
$|\<\Omega\<|$ is an infinite regular cardinal
(e.g., $\aleph_0),$ and if $G$ is a subgroup
of $S=\Sym(\Omega)$ such that every orbit of $\Omega$ under $G$
has cardinality $<|\<\Omega\<|,$ then $G\not\approx_{|\<\Omega\<|}S.$
Given $U\subseteq S$ of cardinality $<|\<\Omega\<|,$
let us define the ``distance'' between elements
$\alpha,\beta\in\Omega$ to be the length of the shortest group word in
the elements of $G\cup U$ that carries $\alpha$ to $\beta,$
or to be~$\infty$ if there is no such word.
It is not hard to see from our assumptions on
the orbits of $G$ and the cardinality of $U$ that for
each $\alpha\in\Omega$ and each positive integer $n,$ there are
$<|\<\Omega\<|$ elements of $\Omega$ within distance $n$ of $\alpha,$
hence this distance function on $\Omega$ has no finite bound.
Now if an element of $\langle\<G\cup U\rangle$ is expressible
as a word of length $n$ in elements of $G\cup U,$
it will move each element of $\Omega$ a distance $\leq n.$
Thus, taking $f\in\Sym(\Omega)$ which moves points by unbounded
distances, we have $f\notin\langle\<G\cup U\rangle;$
so $\langle\<G\cup U\rangle\neq S.$

The next definition gives a name to the concept, used in the
above proof, of a metric under which points may have
distance $\infty,$ and introduces some related terminology
and notation.

\begin{definition}\label{D.metric}
In this definition, $P$ will denote
$\{r\in\mathbb R:r\geq 0\}\cup\{\infty\},$ ordered in the obvious way.

A {\em generalized metric} on a set $\Omega$ will mean a function
$d\<{:}\;\Omega\times\Omega\to P,$
satisfying the usual definition of a metric, except for this
generalization of its value-set.

If $d$ is a generalized metric on $\Omega,$ then
for $r\in P,$ $\alpha\in\Omega,$ we will write
$B_d(\alpha,r)$ for the open ball of radius $r$ about $\alpha,$
$\{\beta\in\Omega:d(\alpha,\beta)<r\}.$
If $\kappa$ is an infinite cardinal, we will call
the generalized metric space $(\Omega,d)$ {\em $\kappa\!\<$-uncrowded}
if for every $\alpha\in\Omega$ and $r<\infty,$
the ball $B_d(\alpha,r)$ has cardinality $<\kappa\<.$
We will call $(\Omega,d)$ {\em uniformly $\kappa\!\<$-uncrowded}
if for every $r<\infty,$ there exists a $\lambda<\kappa$
such that for all $\alpha\in\Omega\<,$ $|B_d(\alpha,r)|\leq\lambda.$
{\rm(}For brevity we will, in these situations, also often refer to
the~{\em metric} $d$ as being $\kappa\!\<$-uncrowded
or uniformly $\kappa\!\<$-uncrowded.{\rm)}

Given two generalized metrics $d$ and $d'$ on $\Omega,$
we shall write $d'\leq d$ if $d'(\alpha,\beta)\linebreak[0]
\leq d(\alpha,\beta)$ for all $\alpha,\beta\in\Omega\<.$

For $g\in\Sym(\Omega)$ and $d$ a generalized metric
on $\Omega,$ we define
$${||g||}_d\ =\ %
\textstyle\sup_{\strt\alpha\in\Omega}\ d(\alpha,\alpha\<g)\;
\in\;P.$$
If ${||g||}_d<\infty,$ then the
permutation $g$ will be called {\em bounded} under $d.$
\end{definition}

The case of these concepts that we will be most concerned with in
subsequent sections is that in which $|\<\Omega\<|=\kappa=\aleph_0.$

Observe that if a group $G\leq\Sym(\Omega)$ has all
orbits of cardinality $<|\<\Omega\<|,$ then giving $\Omega$ the
generalized metric under which distinct points in the same orbit
have distance $1$ and points in different orbits have distance
$\infty,$ we get a $|\<\Omega\<|\!\<$-uncrowded generalized
metric with respect to which $G$ acts by bounded permutations --
this is the $U=\varnothing$ case of the construction sketched in the
second paragraph of this section.
Thus, if we can change the hypothesis of that construction
from ``$G$ has orbits of cardinality $<|\<\Omega\<|$'' to ``$G$
acts by bounded permutations with respect to a
$|\<\Omega\<|\!\<$-uncrowded generalized metric'' (the condition we
deduced held for $\langle\<G\cup U\rangle),$
we will have a stronger statement.
This is done in the next theorem.
If the cardinalities of the orbits of $G$ have a {\em common} bound
$\lambda<|\<\Omega\<|,$ then the above construction gives a
{\em uniformly} $|\<\Omega\<|\!\<$-uncrowded generalized metric.
The theorem will generalize that case as well.

Note that if a generalized metric space $(\Omega,d)$ is
$|\<\Omega\<|\!\<$-uncrowded, the function $d$ is necessarily
unbounded, i.e., takes on values exceeding every positive real number.

\begin{theorem}\label{T.bddG}
Suppose $\Omega$ is an infinite set, $\kappa$ a regular cardinal
$\leq|\<\Omega\<|,$ $d$ a $\kappa\!\<$-uncrowded {\rm(}\!\<respectively
a uniformly $\kappa\!\<$-uncrowded\/{\rm)} generalized metric
on $\Omega,$ and $G$ a subgroup of $S=\Sym(\Omega)$ consisting
of elements bounded with respect to $d.$

Then for any subset $U\subseteq\Sym(\Omega)$ of cardinality
$<\kappa,$ there exists a generalized metric $d'\leq d$ on $\Omega,$
again $\kappa\!\<$-uncrowded {\rm(}\!\<respectively, uniformly
$\kappa\!\<$-uncrowded\/{\rm)}, such that every element
of $U,$ and hence every element of $\langle\<G\cup U\rangle,$
is bounded with respect to $d'.$

Thus, for subgroups $G\leq S,$ the property that there exists a
{\rm(}\!uniformly{\rm)} $\kappa\!\<$-uncrowded
generalized metric on $\Omega$ with respect to which every element
of $G$ acts by bounded permutations is preserved under
passing to groups $H\preq_\kappa G,$ and in particular, under
passing to groups $H\approx_\kappa G.$

\end{theorem}\begin{proof}
Given $\Omega,$ $\kappa,$ $d$ as in the first sentence and $U$ as in
the second, let us define
$d'(\alpha,\beta)$ for $\alpha,\beta\in\Omega,$
to be the infimum, over all finite sequences of the form

\begin{xlist}\item\label{x.a0-an}
$(\alpha_0,\;\alpha_1,\,\ldots\,,\,\alpha_{2n+1})$~ where
$\alpha=\alpha_0,$ $\alpha_{2n+1}=\beta,$ and where for $i$ odd,
$\alpha_{i+1}\in\alpha_i\<(U\cup U^{-1}),$
\end{xlist}
of the sum

\begin{xlist}\item\label{x.d(a,b)}
$d(\alpha_0,\alpha_1)\,+\,1\,+\,d(\alpha_2,\alpha_3)\,+\,1\,+\,
\ldots\,+\,1\,+\,d(\alpha_{2n},\alpha_{2n+1}).$
\end{xlist}
This infimum is $\leq d(\alpha,\beta)$ because the set of sequences
over which it is taken includes the sequence $(\alpha,\beta).$
We also see from~(\ref{x.d(a,b)}) that whenever
$d'(\alpha,\beta)\neq d(\alpha,\beta),$
we have $d'(\alpha,\beta)\geq 1,$ hence $d'(\alpha,\beta)$
is nonzero for $\alpha\neq\beta.$
Symmetry and the triangle inequality are immediate from the
definition, and each $g\in U$ satisfies ${||g||}_{d'}\leq 1$
because elements $\alpha$ and $\alpha\<g$ are connected
by the sequence $(\alpha,\alpha,\alpha\<g,\alpha\<g).$
Thus, $d'$ is a generalized metric $\leq d,$ with respect to which
every element of $U$ is bounded; it remains to show that $d'$
is again (uniformly) $\kappa\!\<$-uncrowded.

So consider a ball of finite radius, $B_{d'}(\alpha,r).$
An element $\beta$ lies in this ball if and only if there is a
sequence~(\ref{x.a0-an}) for which the sum~(\ref{x.d(a,b)}) is $<r.$
But~(\ref{x.d(a,b)}) has $n$ summands equal to $1,$ so given $r<\infty,$
there are only finitely many values of $n$ that need to be
considered; so it suffices to show that for fixed $n$ and $\alpha,$
the number of sequences~(\ref{x.a0-an}) making (\ref{x.d(a,b)})
less than $r$ is $<\kappa,$
and in the uniform case has a bound $<\kappa$ depending only on $r.$
Now for $0\leq i<2n+1,$ if we are given $\alpha_i,$ then the number of
possibilities for $\alpha_{i+1}$ consistent with (\ref{x.d(a,b)}) being
less than $r$ is $\leq|B_d(\alpha_i,r)|<\kappa$
if $i$ is even, while it is $\leq|\,U\cup U^{-1}|<\kappa$ if $i$ is odd.
By the regularity of $\kappa,$ it follows that
the number of possible sequences~(\ref{x.a0-an}) of
length $2n+1$ making (\ref{x.d(a,b)})
less than $r$ and starting with a given $\alpha_0$ is $<\kappa\<.$
Moreover, if we have a bound on $|B_d(\alpha_i,r)|$ independent
of $\alpha_i,$ we also clearly get a bound on the above cardinal
independent of $\alpha,$ as required.

The final assertions, concerning the relations
$\preq_\kappa$ and $\approx_\kappa,$ clearly follow.
\end{proof}

{\it Remark.} What if in the above theorem we weaken the assumption
that $U$ has cardinality $<\kappa$ to say that it
has cardinality $\leq\kappa$?

If $\kappa=\aleph_0,$ we get exactly the same conclusions, since the
result of~\cite{FG} cited in the proof Lemma~\ref{L._0<=>_1}(i)
lets us replace any countable $U$ by a set of cardinality~$2.$
Can we see this stronger assertion without calling on~\cite{FG}?
Yes, by a slight modification of the proof of our theorem:
We write the countable set $U$ as $\{f_1, f_2,\<\ldots\<\},$ and
replace each of the ``$\!1\!$''s in~(\ref{x.d(a,b)}) by a value
$N$ such that $\alpha_{i+1}=\alpha_i\<f_N^{\pm 1}.$
Thus in the final step of the proof, the number of
choices of $\alpha_{i+1}$ that can follow
$\alpha_i$ for $i$ odd will still satisfy a bound
below $\aleph_0,$ namely~$2r.$

For an uncountable regular cardinal $\kappa,$ the sort of generalized
metric we have introduced is not really the best tool.
(Indeed, if $(\Omega,d)$ is $\kappa\!\<$-uncrowded, then for every
$\alpha\in\Omega,$ ``most'' elements of $\Omega$ must be at distance
$\infty$ from $\alpha,$ and the main import of $d$ lies in the
equivalence relation of having distance $<\infty;$ so such a metric
is hardly a significant generalization of an equivalence relation.)
What is more useful then is the concept of a
$\kappa\cup\{\infty\}\!\<$-valued ultrametric, where the symbol $\infty$
is again taken as greater than all other values of the metric.
Defining in the obvious way what it means for such an ultrametric
to be (uniformly) $\kappa\!\<$-uncrowded, one can apply the same method
as above when $|U|=\kappa$, with the
summation~(\ref{x.d(a,b)}) replaced by a supremum.
Since we will not be looking at this situation, we leave
the details to the interested reader.

The above theorem, in generalizing the argument sketched at the
beginning of this section, discarded the explicit connection with
cardinalities of orbits.
The next result records that connection.

Let us understand a {\em partition} of a set $\Omega$ to mean a set
$A$ of disjoint nonempty subsets of $\Omega$ having $\Omega$ as union.
If $A$ is a partition of $\Omega$ and $S=\Sym(\Omega),$ we define

\begin{xlist}\item\label{x.def_SA}
$S_{(A)}\;=\;\{\<f\in S\,:\,
(\<\forall\,\Sigma\in A)\;\;\Sigma f=\Sigma\,\}.$
\end{xlist}
(This is an extension of the notation $S_{(\Sigma)}$ recalled
in the preceding section.)

\begin{theorem}\label{T.no_uncrwdd}
Let $\Omega$ be an infinite set, $A$ a partition of $\Omega,$
$S=\Sym(\Omega),$ $G=S_{(A)},$ and let $\kappa$
be an infinite regular cardinal $\leq|\<\Omega\<|.$
Then
\\[2pt]
{\rm(a)}~ If some member of $A$ has cardinality $\geq\kappa,$
then there is no $\kappa\!\<$-uncrowded generalized metric
on $\Omega$ with respect to which all members of $G$ are bounded.
\\[2pt]
{\rm(b)}~ If all members of $A$ have cardinalities $<\kappa,$ but
there is no common bound $\lambda<\kappa$ for those cardinalities,
then there is a $\kappa\!\<$-uncrowded generalized
metric with respect to which all elements of $G$ are bounded, but no
uniformly $\kappa\!\<$-uncrowded generalized metric with this property.
\\[2pt]
{\rm(c)}~ If all members of $A$ have cardinalities $\leq\lambda$ for
some $\lambda<\kappa,$ then there is a uniformly $\kappa\!\<$-uncrowded
generalized metric with respect to which all elements of $G$ are
bounded.
\vspace{2pt}

Thus, by the last sentence
of Theorem~\ref{T.bddG}, for partitions $A,\,B$ of $\Omega$
falling under distinct cases above, we have
$S_{(A)}\not\approx_\kappa S_{(B)}.$
More precisely, if $A$ falls under a later case than $B,$
then $S_{(A)}\not\sucq_\kappa S_{(B)}.$
\end{theorem}\begin{proof}
To show~(a), let $\Sigma\in A$ have cardinality $\geq\kappa$ and
let $d$ be any $\kappa\!\<$-uncrowded generalized metric on $\Omega\<.$
To construct an element of $G$ which is unbounded with respect
to $d,$ let us choose elements
$\alpha_j,\beta_j\in\Sigma$ for each positive integer $j$ as follows:
Assuming the elements with subscripts $i<j$ have been chosen, take for
$\alpha_j$ any element of $\Sigma$ distinct from all of these.
Since $B_d(\alpha_j,j)$ has cardinality $<\kappa\leq|\<\Sigma\<|,$
the set $\Sigma-B_d(\alpha_j,j)-
\{\alpha_1,\ldots,\alpha_{j{-}1},\beta_1,\ldots,\beta_{j{-}1}\}$
is nonempty; let $\beta_j$ be any element thereof.
Once all $\alpha_j$ and $\beta_j$ are chosen, let $f\in G=S_{(A)}$
interchange $\alpha_j$ and $\beta_j$ for
all $j,$ and fix all other elements of $\Omega\<.$
Since $d(\alpha_j,\beta_j)\geq j$ for each $j,$ $f$ is unbounded
with respect to $d.$

The final assertion of~(b) is shown similarly:
If $d$ is uniformly $\kappa\!\<$-uncrowded, then
for each positive integer $n$ we can find
a cardinal $\lambda_n<\kappa$ such that all balls of
radius $n$ contain $\leq\lambda_n$ elements.
On the other hand, the assumption on $A$ allows us to choose for each
$n$ a set $\Sigma_n\in A$ with more than $\lambda_n+(2n-2)$ elements.
Assuming $\alpha_1,\ldots,\alpha_{n-1},\beta_1,\ldots,\beta_{n-1}$
have been chosen, we take for $\alpha_n$ any element of
$\Sigma_n-\{\alpha_1,\ldots,\alpha_{n-1},\linebreak[0]
\beta_1,\ldots,\beta_{n-1}\}$ and for
$\beta_n$ any element of $\Sigma_n-B_d(\alpha_n,n)-\linebreak[-3]
\{\alpha_1,\ldots,\alpha_{n-1},\beta_1,\ldots,\beta_{n-1}\},$
and finish the argument as before.

To get the positive assertions of~(b) and (c), we define a generalized
metric $d_A$ on $\Omega$ by
letting $d_A(\alpha,\beta)=1$ if $\alpha$ and
$\beta$ are in the same member of $A,$ and $\infty$ otherwise.
This is clearly $\kappa\!\<$-uncrowded,
respectively uniformly $\kappa\!\<$-uncrowded,
and all elements of $G$ are bounded by $1$ under $d_A.$

The conclusions of the final paragraph are straightforward.
\end{proof}

The three cases of the above theorem will be used to separate
the first three of the four situations described in the abstract.
One may ask whether the remaining case can be treated similarly.
For parallelism, one might call a generalized
metric ``absolutely uncrowded'' if all balls of finite
radius are singletons, i.e., if the distance between any
two distinct points is $\infty,$ and then note
that the trivial group is the unique
group of permutations whose elements are bounded
with respect to the absolutely uncrowded generalized metric.
However, the property of acting by bounded permutations
with respect to the unique absolutely uncrowded metric is certainly
not preserved under adjunction of finitely many elements, i.e.,
is not an $\approx_{\aleph_0}\!$-invariant.
Rather than any result of this sort, the property of countability will
separate this fourth equivalence class from the others.

\section{The function topology.}\label{S.func_top}

If $\Omega$ is an infinite set and we regard it as a discrete
topological space, then the set $\Omega^\Omega$ of all functions
$\Omega\to\Omega$
becomes a topological space under the function topology.
In this topology, a subbasis of open sets is given by the sets
$\{f\in\Omega^\Omega:\alpha f=\beta\}$ $(\alpha,\beta\in\Omega\<).$
The closure of a set $U\subseteq\Omega^\Omega$
consists of all maps $f$ such that, for every finite
subset $\Gamma\subseteq\Omega,$ there exists an element of $U$
agreeing with $f$ at all members of $\Gamma.$
It is immediate that composition of maps is continuous in this topology.

The group $S=\Sym(\Omega)$ is not closed in $\Omega^\Omega$ in the
function topology.
For instance when $\Omega=\omega,$ we see that the sequence
of permutations $(0,1),\ (0,1,2),\ \ldots,\ \linebreak[0]
(0,\ldots\<,n),\ldots$ (cycle notation)
converges to the map $n\mapsto n+1,$ which is not surjective.
Nevertheless, when restricted to $S,$ this topology makes
$(~~)^{-1}$ as well as composition continuous; indeed,
$\{f\in S:\alpha f=\beta\}^{-1}=\{f\in S:\beta f=\alpha\}.$

Given a subset $U\subseteq S,$ we shall write $\cl(U)$ for the closure
of $U$ in $S$ (not in $\Omega^\Omega$!) under the function topology.
The fact that $S$ is not closed in $\Omega^\Omega$ has
the consequence that if one wants to prove the existence
of an element $f\in\cl(U)$ behaving in some desired fashion, one
cannot do this simply by finding elements
of $U$ that show the desired behavior
at more and more elements of $\Omega,$ and saying ``take the limit'';
for the limit may be an element of $\Omega^\Omega$ which is not in $S.$
However, there is a standard way of getting around this
difficulty, ``the method of going back and forth''.
One constructs elements of $U$ which not only agree on
more and more elements of $\Omega,$ but whose inverses also
agree on more and more elements.
Taking the limit, one thus gets a map and also an inverse to that map.
Cf.\ \cite[\S\S9.2, 16.4]{BhMMN} for examples
of this method, and some discussion.
The next result, a formalization of this idea, will be used
at several points below.

\begin{lemma}\label{L.b+f}
Suppose that $\Omega=\{\eps_0,\eps_1,\ldots\,\}$ is a countably
infinite set, and that $g_0,\,g_1,\,\ldots\in S=\Sym(\Omega)$ and
$\Gamma_0,\Gamma_1,\,\ldots\subseteq\Omega$ are such that for all $j>0,$

\begin{xlist}\item\label{x.0-j.1,g-1}
$\{\<\eps_0,\ldots,\eps_{j-1\<}\}\,\cup\,
\{\<\eps_0\,g_{j-1}^{-1},\ldots,\eps_{j-1}\,g_{j-1}^{-1\<}\}\;
\subseteq\;\Gamma_j,$
\end{xlist}
and

\begin{xlist}\item\label{x.gjinS_*Gg}
$g_j\in S_{(\Gamma_j)}\;g_{j-1}.$
\end{xlist}
Then the sequence $(g_j)_{j=0,1,\ldots}$ converges in $S.$
\end{lemma}\begin{proof}
Let $i\geq 0.$

For all $j>i,$ the conditions $\eps_i\in\Gamma_j$
and~(\ref{x.gjinS_*Gg}) imply that $\eps_i\,g_j=\eps_i\,g_{j-1}.$
Thus, the sequence $(g_j)$ is eventually constant on $\eps_i.$

Likewise,~(\ref{x.gjinS_*Gg}) and the condition
$\eps_i\,g_{j-1}^{-1}\in\Gamma_j$ imply that
$\eps_i\,g_j^{-1}=\eps_i\,g_{j-1}^{-1};$ hence the element
of $\Omega$ carried to $\eps_i$ by $g_j$ is the same for all $j>i.$

Since the first conclusion holds for all $\eps_i\in\Omega,$
the sequence $(g_j)_{j=0,1,\ldots}$ converges to an element
of $\Omega^\Omega,$ which is one-to-one because all the $g_j$ are.
Applying the second conclusion, we see that each
$\eps_i$ is in the range of $g,$ so $g\in\Sym(\Omega).$
\end{proof}

We note some elementary
facts about closures of subgroups in the function topology.

\begin{lemma}\label{L.clG_*G}
Suppose $\Omega$ is a set and $G$ a subgroup of $S=\Sym(\Omega).$
Then
\\[2pt]
{\rm(i)}~ $\cl(G)$ is also a subgroup of $S.$
\\[2pt]
{\rm(ii)}~ $G$ and $\cl(G)$ have the same orbits in $\Omega\<.$
\\[2pt]
{\rm(iii)}~ If $\Gamma$ is a finite subset of $\Omega,$ then
$\cl(G)_{(\Gamma)}=\cl(G_{(\Gamma)}).$
\end{lemma}\begin{proof}
Statement~(i) is an immediate consequence of the continuity
of the group operations.

 From the characterization of the closure of a
set in our topology, we see that for $\alpha,\,\beta\in\Omega,$
the set $\cl(G)$ will contain elements carrying $\alpha$
to $\beta$ if and only if $G$ does, from which~(ii) is clear.

The direction $\supseteq$ in~(iii) follows by applying~(ii) to the
orbits of elements of $\Gamma.$
(Finiteness of $\Gamma$ is not needed for this direction.)
To get $\subseteq,$ assume $f\in\cl(G)_{(\Gamma)}.$
Since $f\in\cl(G),$ every neighborhood of $f$ contains elements of $G.$
But as $f$ fixes all points of the finite set $\Gamma,$ every
sufficiently small neighborhood of $f$ consists of elements which do
the same, hence every such neighborhood contains points
of $G_{(\Gamma)};$ so $f\in\cl(G_{(\Gamma)}).$
\end{proof}

The above lemma has the consequence that once we show
(for $\Omega$ countable) that
the $\approx_{\aleph_0}\!$-class of a closed subgroup of $\Sym(\Omega)$
is determined by which of conditions~(i)--(iv) in our abstract
hold, we can also say for an arbitrary subgroup $G\leq\Sym(\Omega)$
that the $\approx_{\aleph_0}\!$-class of $\cl(G)$ is determined
in the same way by which of those conditions $G$ satisfies.

The subgroups of $\Sym(\Omega)$ closed in the function topology are
known to be precisely the automorphism groups of the finitary relational
structures on $\Omega\<.$
(Indeed, one may take the $n\!$-ary relations in such a
structure, for each $n,$ to be all orbits of
$n\!\<$-tuples of elements of $\Omega$ under the group.)
But we shall not make use of this fact here.

(Incidentally, $\Sym(\Omega)$ is also not {\em open} in $\Omega^\Omega.$
It is easy to give a sequence of non-injective or non-surjective
maps in which the failures of injectivity or surjectivity ``drift off
to infinity'', so that the limit is a bijection, e.g., the identity.)

\section{Infinite orbits.}\label{S.inf_orbs}

In this and the next three sections (and with minor exceptions,
in subsequent sections as well), we shall restrict attention to
the case of countable $\Omega\<.$
When an enumeration of its elements is required, we shall write
\begin{xlist}\item\label{x.enumerate}
$\Omega\ =\ \{\<\eps_i\;:\;i\in\omega\<\}.$
\end{xlist}
References to limits etc.\ in $S=\Sym(\Omega)$ will always refer to the
function topology; in particular, a closed subgroup
of $S$ will always mean one closed in $S$ under that topology.
The symbols $\preq$ and $\approx$ will mean $\preq_{\strt\aleph_0,S}$
and $\approx_{\strt\aleph_0,S}$ respectively.

We shall show in this section that if $G$ is a closed subgroup
of $S$ such that

\begin{xlist}\item\label{x.inf_orbs}
For every finite subset $\Gamma\subseteq\Omega,$ the subgroup
$G_{(\Gamma)}$ has at least one infinite orbit in $\Omega,$
\end{xlist}
then $G\approx S.$
Our proof will make use of the following result of
Macpherson and Neumann:

\begin{xlist}\item\label{x.M&N.full}
\cite[Lemma~2.4]{DM&PN} (cf.\ \cite[Lemma~3]{Sym_Omega:1}):~
Suppose $\Omega$ is an infinite set and $H$ a subgroup of
$\Sym(\Omega),$ and suppose there exists a
subset $\Sigma\subseteq\Omega$ of the same cardinality as $\Omega,$
such that $H_{\{\Sigma\}}$ (i.e., $\{\<f\in H:\Sigma f=\Sigma\<\})$
induces, under restriction to
$\Sigma,$ the full permutation group of $\Sigma\<.$
Then there exists $x\in\Sym(\Omega)$ such that
$\Sym(\Omega)=\langle H\cup\{x\}\<\rangle.$
\end{xlist}
(This is stated in~\cite{DM&PN} and~\cite{Sym_Omega:1} for
the case where $\Sigma$ is a {\em moiety,} i.e., a set of
cardinality $|\<\Omega\<|$ such that $\Omega-\Sigma$ also has
cardinality $|\<\Omega\<|\<.$
But if the hypothesis of~(\ref{x.M&N.full}) holds for some
$\Sigma $ of cardinality $|\<\Omega\<|,$ it clearly also holds
for a subset of $\Sigma$ which is a moiety,
so we may restate the result as above.)

We will also use the following fact.
We suspect it is known, and would appreciate learning of any reference.
(A similar technique, but not this result,
occurs in \cite{Such84} and \cite{Such85}.)

\begin{lemma}\label{L.loc_bdd}
Let us call a permutation $g$ of the set $\omega$ of natural
numbers {\em local} if for every $i\in\omega$ there exists $j>i$
in $\omega$ such that $g$ carries $\{0,\ldots,j\<{-}\<1\}$ to itself.

Then every permutation $f$ of $\omega$ is a product
$gh$ of two local permutations.
\end{lemma}\begin{proof}
Given $f\in\Sym(\omega),$ let us choose integers
$0=a(0)<a(1)<a(2)<\ldots$ recursively, by letting each $a(i)$ be any
value $>a(i{-}1)$ such that $\{0,\ldots,a(i{-}1){-}1\}\,f\linebreak[0]
\,\cup\,\{0,\ldots,a(i{-}1)\<{-}\<1\}\,f^{-1}
\subseteq\{0,\ldots,a(i)\<{-}\<1\}.$
Let $\Sigma_{-1}=\varnothing,$ and for $i\geq 0$ let
$\Sigma_i=\{a(i),\,a(i){+}1,\,\ldots,\,a(i{+}1){-}1\}.$
Thus the set $A=\{\Sigma_i:i\geq 0\}$ is a
partition of $\omega$ into finite subsets, such that for each $i\geq 0$
one has $\Sigma_i f\subseteq\Sigma_{i-1}\cup\Sigma_i\cup\Sigma_{i+1}.$
Note that for each $i\geq 0,$ the number of elements
which are moved by $f$ from $\Sigma_{i-1}$ into $\Sigma_i$ is equal to
the number that are moved from $\Sigma_i$ into $\Sigma_{i-1}$
(since these are the elements of $\omega$
that are moved ``past $a(i)-1/2$'' in the upward, respectively
the downward direction).

We shall now construct a permutation $g$ such that $g$ carries each set
$\Sigma_{2i}\<\cup\<\Sigma_{2i+1}$ $(i\geq 0)$ into itself, and
$g^{-1} f$ carries each set $\Sigma_{2i-1}\<\cup\<\Sigma_{2i}$
$(i\geq 0)$ into itself; thus each of these permutations will
be local, and they will have product $f,$ as required.
To do this let us, for each $i,$ pair elements $\alpha$ that $f$ carries
from $\Sigma_{2i}$ upward into $\Sigma_{2i+1}$ with elements $\beta$
that it carries from $\Sigma_{2i+1}$ downward into $\Sigma_{2i}$
(having seen that the numbers of such elements are equal), and let $g$
exchange the members of each such pair, while fixing other elements.
Clearly $g$ preserves the sets $\Sigma_{2i}\<\cup\<\Sigma_{2i+1}.$
It is not hard to verify that if we now look at $\omega$ as
partitioned the other way, into the intervals
$\Sigma_{2i-1}\<\cup\<\Sigma_{2i},$ then the $g$
we have constructed has the property that for every $\alpha\in\omega,$
the element $\alpha\<g$ lies in the same interval
$\Sigma_{2i-1}\<\cup\<\Sigma_{2i}$ as does $\alpha f.$
Hence $g^{-1}f$ preserves each interval
$\Sigma_{2i-1}\<\cup\<\Sigma_{2i},$ completing the proof.
\end{proof}

We shall now prove a generalization of~(\ref{x.M&N.full}), assuming
$\Omega$ countable.
To motivate the statement, note that in the countable
case of~(\ref{x.M&N.full}), if we enumerate the elements
of $\Sigma$ as $\alpha_0,\<\alpha_1,\<\ldots\<,$
then the hypothesis implies that we can choose elements $g\in H$
in ways that allow us infinitely many choices for
$\alpha_0\<g,$ for each such choice infinitely many choices
for $\alpha_1\<g,$ etc..
But the hypothesis of~(\ref{x.M&N.full})
is much stronger than this, since it specifies
that the set of choices for $\alpha_0 g$ include all the $\alpha_i,$
that the choices for $\alpha_1\<g$ then include all
$\alpha_i$ other than the one chosen to be $\alpha_0\<g,$ etc..
The next result says that we can get the same
conclusion without such a strong form of the hypothesis.

\begin{lemma}\label{L.D_i}
Let $\Omega$ be a countably infinite set and $G$ a subgroup of
$\Sym(\Omega),$ and suppose there exist a sequence
$(\alpha_i)_{i\in\omega}\in\Omega^\omega$ of distinct elements,
and a sequence of
nonempty subsets $D_i\subseteq\Omega^i$ $(i\in\omega),$ such that
\\[2pt]
{\rm(i)}~ For each $i\in\omega$ and
$(\beta_0,\ldots,\beta_{i})\in D_{i+1},$ we have
$(\beta_0,\ldots,\beta_{i-1})\in D_i;$
\\[2pt]
{\rm(ii)}~ For each $i\in\omega$ and
$(\beta_0,\ldots,\beta_{i-1})\in D_i,$ there exist infinitely
many elements $\beta\in\Omega$ such that
$(\beta_0,\ldots,\beta_{i-1},\beta)\in D_{i+1};$ and
\\[2pt]
{\rm(iii)}~ If $(\beta_i)_{i\in\omega}\in\Omega^\omega$ has the property
that $(\beta_0,\ldots,\beta_{i-1})\in D_i$ for each $i\in\omega,$
then there exists $g\in G$ such that $(\beta_i)=(\alpha_i\,g)$
in $\Omega^\omega.$
\vspace{2pt}

Then $G\approx S.$
\end{lemma}\begin{proof}
Let us note first that our hypotheses imply that for
$(\beta_0,\ldots,\beta_{i-1})\in D_i,$ the entries $\beta_j$ are all
distinct.
For from~(i) and~(ii) we see that such an $i\!\<$-tuple can be extended
to an $\omega\!\<$-tuple as in~(iii), and by~(iii)
this $\omega\!\<$-tuple
is the image under a group element of the $\omega\!\<$-tuple
of distinct elements $(\alpha_i).$

We shall now construct recursively, for $i=0,1,\ldots\,,$
finite sets $E_i\subseteq D_i.$
For each $i,$ the elements of $E_i$ will be denoted
$e(n_0,\ldots,n_r;\,\pi_1,\ldots,\pi_r)$ with one such element
for each choice of a sequence of natural numbers
$0=n_0<n_1<\ldots<n_r=i$ and a sequence of permutations
$\pi_m\in\Sym(\{n_{m-1},\ldots,n_m{-}1\})$ $(1\leq m\leq r).$
(Note that each $e(n_0,\ldots,n_r;\,\pi_1,\ldots,\pi_r),$
since it belongs to $D_i,$ is an
$i\!\<$-tuple of elements of $\Omega,$ where $i=n_r;$
but we shall not often write it explicitly as a string of elements.
Nevertheless, we shall refer to the $i$ elements of $\Omega$ comprising
this $i$\!-tuple as its {\em components}.)

We start the recursion with $E_0=D_0,$ which is necessarily
the singleton consisting of the unique length-$\!\<0$ sequence.
Assuming $E_0,\ldots,E_{i-1}$ given, we choose an arbitrary order in
which the finitely many $i\!\<$-tuples in $E_i$ are to be chosen.
When it comes time to choose the $i\!\<$-tuple
$e(n_0,\ldots,n_r;\,\pi_1,\ldots,\pi_r)\in E_i,$ we define its
initial substring of length $n_{r-1}$ to be the $n_{r-1}\!$-tuple
$e(n_0,\ldots,n_{r-1};\,\pi_1,\ldots,\pi_{r-1})\in E_{n_{r-1}}.$
We then extend this to an element of $D_{n_r}$ in any way such that
its remaining $n_r-n_{r-1}$ components are distinct from all components
of all elements of $E_0\cup\ldots\cup E_{i-1},$ and from all components
of those elements of $E_i$ that have been chosen so far.
This is possible by $n_r-n_{r-1}$ applications of condition~(ii) above:
at each step, when we extend a member of a set $D_j$ to a member of
the next set $D_{j+1}$ $(n_{r-1}\leq j<n_r)$ we have infinitely many
choices available for the last component,
and only finitely many elements to avoid.

Once the sets $E_i$ are chosen for all $i,$ let us define an element
$s\in\Sym(\Omega)$ to permute, in the following way, those elements
of $\Omega$ that occur as components in the members of $\bigcup_i E_i.$
(On the complementary subset of $\Omega$ we let $s$ behave
in any manner, e.g., as the identity.)

\begin{xlist}\item\label{x.*b_j.s}
For each $(\beta_j)_{0\leq j<n_r}=
e(n_0,\ldots,n_r;\,\pi_1,\ldots,\pi_r)\in E_i,$
we let $s$ act on its last $n_r-n_{r-1}$
components, $\beta_{n_{r-1}}, \beta_{n_{r-1}+1},\ldots,\beta_{n_r-1},$
by\linebreak[4] $\beta_j\<s =\nolinebreak\beta_{j\pi_r}.$
\end{xlist}
That is, we let $s$ permute the
elements $\beta_{n_{r-1}},\ldots,\beta_{n_r-1}$ by
``acting as $\pi_r$ on their subscripts''.
Note that (by the choices made in the last paragraph),
for each $j\in\{n_{r-1},\ldots,n_r{-}1\},$
the occurrence of $\beta_j$ as a component of
$e(n_0,\ldots,n_r;\,\pi_1,\ldots,\linebreak[0]\pi_r)$ is its
first appearance
among the components of the elements we have constructed, and that it
is distinct from the elements first appearing as components
of other tuples $e(n'_0,\ldots,n'_r;\,\pi'_0,\ldots,\pi'_r),$
or in other positions of $e(n_0,\ldots,n_r;\,\pi_1,\ldots,\pi_r).$
Thus~(\ref{x.*b_j.s}) uniquely defines $s$ on this set of elements.

Consider now any permutation of $\{\alpha_i\}$ of the form
$\alpha_i\mapsto\alpha_{i\pi}$ where $\pi$ is a {\em local} permutation
of $\omega$ (in the sense of Lemma~\ref{L.loc_bdd}).
We claim that there exists $g\in G$ such that the element $s$
constructed above ``acts as $\pi$ on the subscripts'' of the
image sequence
$(\alpha_i\,g)_{\strt i\geq 0},$ i.e., such that for all $i\geq 0,$

\begin{xlist}\item\label{x.*a_i.g.s}
$\alpha_i\,g\,s\;=\;\alpha_{i\pi}\,g.$
\end{xlist}

To show this, note that since $\pi$ is local, we can find natural
numbers $0=n_0<n_1<\ldots$ such that $\pi$ carries each of the intervals
$\{n_{m-1},n_{m-1}{+}1,\ldots,n_m{-}1\}$ into itself.
Let us denote the restrictions of $\pi$ to these intervals
by $\pi_m\in\Sym(\{n_{m-1},\linebreak[0]
n_{m-1}{+}1,\ldots,n_m{-}1\})$ $(m\geq 1),$
and consider the tuples

\begin{xlist}\item\label{x.e_...}
$e(n_0)\in E_0,\quad e(n_0,n_1;\,\pi_1)\in E_{n_1},\quad \ldots\;,$\\
\indent \indent \indent \indent
$e(n_0,\ldots,n_m;\,\pi_1,\ldots,\pi_m)\in
E_{n_m},\quad\ldots\;.$
\end{xlist}
Each of these tuples extends the preceding, so there is a sequence
$(\beta_i)\in\Omega^\omega$ of which these tuples are all truncations.
From~(\ref{x.*b_j.s}) we see that the sequence $(\beta_i)$ will
satisfy $\beta_i\<s=\beta_{i\pi}$ for all $i\in\omega.$
Also, by our hypothesis~(iii) and the
condition $E_i\subseteq D_i,$ there exists
$g\in G$ such that $\beta_i=\alpha_i\<g$ for all $i.$
Substituting this into the relation $\beta_i\<s=\beta_{i\pi},$
we get~(\ref{x.*a_i.g.s}), as claimed.

Now~(\ref{x.*a_i.g.s}) can be rewritten as
saying that $g\,s\,g^{-1}$ acts on $\{\alpha_i\}$ by the
map $\alpha_i\mapsto\alpha_{i\pi}.$
In view of Lemma~\ref{L.loc_bdd}, every permutation
of $\{\alpha_i\}$ can be realized as the restriction to that
set of a product of two such permutations,
hence as $g\,s\,g^{-1}h\,s\,h^{-1}$ for some $g,h\in G.$
Thus, the group $H=\langle\<G\cup\{s\}\rangle$ satisfies the hypothesis
of~(\ref{x.M&N.full}) with $\Sigma=\{\alpha_i\}.$
Hence by~(\ref{x.M&N.full}) there exists $x\in\Sym(\Omega)$
such that $\langle\<G\cup\{s,x\}\rangle=S,$
completing the proof of the lemma.
\end{proof}

We now consider a
subgroup $G\leq\Sym(\Omega)$ satisfying~(\ref{x.inf_orbs}).
We shall show how to construct elements $\alpha_i\in\Omega$ and families
$D_i\subseteq\Omega^i$ satisfying conditions~(i) and~(ii) of the above
lemma, and such that if $G$ is closed, condition~(iii) thereof
also holds, allowing us to apply that lemma.

We begin with another recursion, in which we will construct for each
$j\geq 0$ an element $\alpha_j,$ and a finite subset $K_j$ of
$G,$ indexed

\begin{xlist}\item\label{x.g()}
$K_j=\{g(k_0,k_1,\ldots,k_{r-1}):k_0,k_1,\ldots,
k_{r-1},r\in\omega,\,r+k_0+\ldots+\nolinebreak
k_{r-1}=\nolinebreak j\}.$
\end{xlist}

To describe the recursion, assume inductively that $\alpha_i$
and $K_i$ have been defined for all nonnegative $i<j,$
and let $\Gamma_j\subseteq\Omega$ denote the finite set consisting of
the images of $\eps_0,\ldots,\eps_{j-1}$
(cf.~(\ref{x.enumerate})) and of
$\alpha_0,\ldots,\alpha_{j-1}$ under the inverses
of all elements of $K_0\cup\ldots\cup K_{j-1}.$
Let $\alpha_j$ be any element of $\Omega$ having infinite orbit
under $G_{(\Gamma_j)}$ (cf.~(\ref{x.inf_orbs})).
In choosing the elements $g(k_0,k_1,\ldots,k_{r-1})$
comprising $K_j,$ we consider two cases.

If $j=0,$ we have only one element, $g(\<),$ to choose, and we
take this to be the identity element $1\in G.$
A consequence of this choice is that for all larger $j,$
we have $1\in K_0\cup\ldots\cup K_{j-1},$
hence the definition of $\Gamma_j$ above guarantees
that $\eps_0,\ldots,\eps_{j-1}$ and
$\alpha_0,\ldots,\alpha_{j-1}$ themselves lie in $\Gamma_j.$

If $j>0,$ we fix arbitrarily an order in which the
elements of $K_j$ are to be constructed.
When it is time to construct $g(k_0,k_1,\ldots,k_{r-1}),$
let us write $g'=g(k_0,k_1,\ldots,k_{r-2}),$ noting
that this is a member of $K_{j-k_{r-1}-1},$ hence already defined.
We will take for $g(k_0,k_1,\ldots,k_{r-1})$ the result of
left-multiplying $g'$ by a certain element $h\in G_{(\Gamma_{r-1})}.$
Note that whatever value in this group we choose for $h,$
the images of $\alpha_0,\ldots,\alpha_{r-2}$ under $h\<g'$ will
be the same as their images under $g',$ since elements of
$G_{(\Gamma_{r-1})}$ fix $\alpha_0,\ldots,\alpha_{r-2}\in\Gamma_{r-1}.$
On the other hand, we may choose $h$ so that
the image of $\alpha_{r-1}$ under $h\<g'$
is distinct from the images of $\alpha_{r-1}$ under
the finitely many elements of $K_0\cup\ldots\cup K_{j-1},$ and also
under the elements of $K_j$ that have so far been constructed, since
$\alpha_{r-1}$ has infinite orbit under $G_{(\Gamma_{r-1})},$
and there are only finitely many elements that have to be avoided.
So let $g(k_0,k_1,\ldots,k_{r-1})=h\<g'$ be so chosen.

In this way we successively construct the elements of each set $K_j.$
Note that this gives us group elements $g(k_0,\ldots,k_{i-1}\<)$ for all
$i,\ k_0,\ldots,k_{i-1}\in\omega.$
We can thus define, for each $i\in\omega,$

\begin{xlist}\item\label{x.def_D_i}
$D_i=\{(\alpha_0\<g,\ldots,\alpha_{i-1}\<g):
g=g(k_0,\ldots,k_{i-1})$ for some $k_0,\ldots,k_{i-1}
\in\nolinebreak\omega\}.$
\end{xlist}

By construction, $g(k_0,\ldots,k_i)$ agrees with $g(k_0,\ldots,k_{i-1})$
on $\alpha_0,\ldots,\alpha_{i-1},$ so chopping off the
last component of an element of $D_{i+1}$ gives an element of $D_i,$
establishing condition~(i) of Lemma~\ref{L.D_i}.
Moreover, any two elements of the form $g(k_0,\ldots,k_i)\in D_{i+1}$
with indices $k_0,\ldots,k_{i-1}$
the same but different last indices $k_i$ act differently
on $\alpha_i,$ so the sets $D_i$ satisfy condition~(ii) of that lemma.
Suppose, now, that $(\beta_i)\in\Omega^\omega$ has the property that
for every $i$ the sequence $(\beta_0,\ldots,\beta_{i-1})$ is in $D_i.$
We see inductively that successive strings
$(\<),\;(\beta_0),\;\dots,\,(\beta_0,\ldots,\beta_{i-1}),\;\ldots$
must arise from unique elements of the forms
$g(\<),\;g(k_0),\;\ldots\,,\;g(k_0,\ldots,k_{i-1}),\;\ldots\,.$
Moreover, by construction each of these group elements
$g(k_0,\ldots,k_i)$ is obtained from the preceding element
$g(k_0,\ldots,k_{i-1})$ by left multiplication by an
element of $G_{(\Gamma_i)},$ where $\Gamma_i$ contains the elements
$\eps_0,\ldots,\eps_{i-1}$ and their inverse images under all the
preceding group elements.
It follows by Lemma~\ref{L.b+f} that if $G$ is closed, the above
sequence converges to an element $g\in G$ whose behavior on
$(\alpha_i)$ is the limit of the behaviors of these elements,
i.e., which sends $(\alpha_i)$ to $(\beta_i),$ establishing
condition~(iii) of Lemma~\ref{L.D_i}.
Hence that lemma tells us that $G\approx S.$

We can now easily obtain

\begin{theorem}\label{T.inf_orbs}
Let $\Omega$ be a countably infinite set, and
$G$ a closed subgroup of $S=\Sym(\Omega).$
Then $G\approx S$ {\rm(}i.e., $S$ is finitely generated over
$G\!\<${\rm)} if and only if $G$ satisfies~{\rm(\ref{x.inf_orbs})}.
\end{theorem}\begin{proof}
We have just seen that~(\ref{x.inf_orbs}) implies $G\approx S.$
On the other hand, if~(\ref{x.inf_orbs}) does not hold,
then for some finite $\Gamma\subseteq\Omega,$
$G_{(\Gamma)}$ has only finite orbits.
Letting $A$ denote the set of these orbits,
we have $G_{(\Gamma)}\leq S_{(A)}.$
But $S_{(A)}$ falls under case~(b) or~(c) of Theorem~\ref{T.no_uncrwdd}
(with $\kappa=\aleph_0)$ while
$S$ falls under case~(a), being determined by the
improper partition of $\Omega\<.$
We thus get
\begin{xlist}\item\label{x.clG<S}
$G\;\approx\;G_{(\Gamma)}\;\leq\;S_{(A)}\;\prec\;S,$
\end{xlist}
where the first relation holds by Lemma~\ref{L.G_*G} and
Lemma~\ref{L._0<=>_1}(i), and the final strict inequality by the last
sentence of Theorem~\ref{T.no_uncrwdd}.
Thus $G\not\approx S.$
\end{proof}

{\em Notes on the development of the above theorem:}
In the proof of Lemma~\ref{L.D_i}, and again in the arguments following
that proof, it might at first appear that our hypotheses
that certain subsets of $\Omega$ were infinite (namely,
in the former case, the set of ``next terms'' extending each member
of $D_i,$ and in the latter, at least one orbit of
$G_{(\Gamma)}$ for each finite set $\Gamma)$
could have been replaced by statements that those sets could be taken
to have large enough finite cardinalities, since at each
step, we had to make only finitely many choices from these sets,
and to avoid only finitely many elements of $\Omega\<.$
But closer inspection shows that we dipped into
these sets for additional elements infinitely many times.
In the proof of Lemma~\ref{L.D_i}, this is because
for fixed $n_0,\dots,n_{r-1}$ there are
infinitely many possibilities for $n_r>n_{r-1},$ and for {\em each}
of these, the construction of $E_{n_r}$ requires extending the elements
$e(n_0,\ldots,n_{r-1};\,\pi_1,\ldots,\pi_{r-1})\in D_{n_{r-1}}$ to
elements of $D_{n_{r-1}+1}.$
Likewise, in (\ref{x.g()}), note that $r\leq j,$ and
each value of $r$ comes up for infinitely many $j,$ so that for each
$r$ we must choose, in the long run, elements of $G_{(\Gamma_{r-1})}$
having infinitely many different effects on $\alpha_{r-1}.$

This spreading out of the choices we made from each infinite
set, into infinitely many clumps of finitely many choices each,
was necessary:  If we had made infinitely many choices at one time
from one of our sets, we would have had infinitely many obstructions
to our choices from the next set, and could not have argued that
those choices could be carried out as required.

Could the two very similar recursive constructions just referred
to have been carried out simultaneously?
In an earlier draft of this note they were.
That arrangement was more efficient (if less transparent as to
what was being proved), and could be considered preferable if
one had no interest except in closed subgroups.
However, the present development yields the
intermediate result Lemma~\ref{L.D_i},
which can be used to show the $\approx\!\<$-equivalence to $S$
of many non-closed subgroups for which, so far as we can see,
Theorem~\ref{T.inf_orbs} is of no help.

For example, consider a partition $A$ of $\Omega$ into a countably
infinite family of countably infinite sets $\Sigma_i,$ and let
$G$ be the group of permutations of $\Omega$ that, for each $i,$ carry
$\Sigma_i,$ into itself, and move only finitely many elements
of that set.
If we choose an arbitrary element $(\alpha_i)\in\prod_i\<\Sigma_i,$
and let $D_i=\Sigma_0\times\ldots\times\Sigma_{i-1}$ for each $i,$
then we see easily that the conditions of Lemma~\ref{L.D_i} hold,
hence that $G\approx S.$

(The same argument works for the subgroup of the above $G$ consisting
of those elements $g$ for which there is a bound independent
of $i$ on the number of elements of $\Sigma_i$ moved by $g.)$

\section{Finite orbits of unbounded size.}\label{S.unbdd_orbs}
In this section, we again let $\Omega$ be a countably infinite set,
and will show that all closed subgroups $G\leq S=\Sym(\Omega)$ for which

\begin{xlist}\item\label{x.unbdd_orbs}
There exists a finite subset $\Gamma\subseteq\Omega$ such that all
orbits of $G_{(\Gamma)}$ are finite, but no such $\Gamma$ for
which the cardinalities of these orbits have a common finite bound,
\end{xlist}
are mutually $\approx\!\<$-equivalent.
The approach will parallel that of the preceding section, but
there are some complications.

First, there is not one natural subgroup that represents this
equivalence class, as $S$ represented the class considered
in the previous section.
Instead we will begin by defining a certain natural
family of closed subgroups of $S$ which we
will prove $\approx\!\<$-equivalent to one another.
Second, we do not have a result from the literature to serve in
the role of~(\ref{x.M&N.full}).
So we will prove such a result.
The fact that a finite symmetric group is not its own commutator
subgroup will complicate the latter task.
(Cf.\ the proof of~(\ref{x.M&N.full}) as
\cite[Lemma~3]{Sym_Omega:1}, which uses the fact, due to Ore~\cite{OO},
that every element of an infinite symmetric group is a commutator.)
So we will prepare for that proof by showing
that certain infinite products of finite symmetric groups within
$S$ are $\approx\!\<$-equivalent to
the corresponding products of alternating groups.
\vspace{2pt}

To define our set of representatives of the $\approx\!\<$-equivalence
class of subgroups of $S$ we are interested in, let

\begin{xlist}\item\label{x.def_P}
$\mathcal P\;=\;\{A: A$ is a partition of $\Omega$ into finite subsets,
and there is no common finite bound on the cardinalities of the
members of $A\,\}\<.$
\end{xlist}
For $A\in\mathcal P$ (and $S_{(A)}$ defined
by~(\ref{x.def_SA})), we see that

\begin{xlist}\item\label{x.SAiso}
$S_{(A)}\;\cong\;\prod_{\strut\Sigma\in A}\Sym(\Sigma).$
\end{xlist}
Note that if a partition $A_1$ is the image of a partition $A_2$ under
a permutation $f$ of $\Omega,$ then $S_{(A_2)}$ is the conjugate of
$S_{(A_1)}$ by $f;$ in particular, $S_{(A_1)}\approx S_{(A_2)}.$
Let us now show more; namely, that

\begin{xlist}\item\label{x.SAeqSB}
$S_{(A)}\;\approx\;S_{(B)}$~ for all $A,B\in\mathcal P.$
\end{xlist}

We claim first that given $A, B\in\mathcal P,$ we can find
two elements $f, g\in S$ such that

\begin{xlist}\item\label{x.Df_or_Dg}
For each $\Sigma\in B$ there exists $\Delta\in A$ such that
$\Sigma\subseteq\Delta f$ or $\Sigma\subseteq\Delta\<g.$
\end{xlist}
Indeed, write $B$ as the disjoint union of any two
infinite subsets $B_1$ and $B_2.$
We may construct $f$ by defining it on one member of $A$ after
another, making sure that each member of $B_1$ ends up within
the image $\Delta f$ of some sufficiently large $\Delta\in A.$
We map those members of $A$ or subsets of members of $A$ that are not
used in this process into the infinite set
$\bigcup_{\Sigma\in B_2}\Sigma\<,$ and
we also make sure to include every element of
$\bigcup_{\Sigma\in B_2}\Sigma$ in the range of
$f,$ so that $f$ is indeed a permutation.
We similarly construct $g$ so that every member of $B_2$
is contained in the image $\Delta\<g$ of some $\Delta\in A.$
Condition~(\ref{x.Df_or_Dg}) is thus satisfied.

For such $f$ and $g,$ we claim that

\begin{xlist}\item\label{x.SB_in_prod}
$S_{(B)}\;\subseteq (f^{-1}S_{(A)}f)\<(g^{-1}S_{(A)}\,g).$
\end{xlist}
Indeed, every element of $S_{(B)}$ can be written as the product
of a member of $S_{(B)}$ which moves only elements of
$\bigcup_{\Sigma\in B_1}\Sigma$ and one which moves only elements
of $\bigcup_{\Sigma\in B_2}\Sigma\<;$ and these can be seen to
belong to $f^{-1}S_{(A)}f$ and to $g^{-1}S_{(A)}\,g$ respectively.

Thus $S_{(B)}\leq\langle\<S_{(A)}\cup\{f,g\}\<\rangle,$
so $S_{(B)}\preq S_{(A)}.$
Since this works both ways, we get~(\ref{x.SAeqSB}), as desired.
\vspace{2pt}

We next prepare for the difficulties concerning
alternating groups versus symmetric groups.
Let $A$ be a partition belonging to $\mathcal P$
which contains infinitely many singletons, and whose other members
are all of cardinality at least $4,$ and let $S_{(A)}^\mathrm{\,even}$
denote the subgroup of $S_{(A)}$ which acts by
an even permutation on each member of $A.$
We shall show that

\begin{xlist}\item\label{x.even=all}
$S_{(A)}^\mathrm{\,even}\;\approx\;S_{(A)}.$
\end{xlist}

To do this, let us list the non-singleton members of $A$ as
$\Sigma_0,\,\Sigma_1,\,\ldots,$ and for each $i$ choose in $\Sigma_i$
four distinct elements, which we name
$\alpha_{4i},\,\alpha_{4i+1},\,\alpha_{4i+2},\,\alpha_{4i+3}.$
Let $B$ denote the partition of $\Omega$ (not belonging
to $\mathcal P)$ whose only nonsingleton subsets are the
two-element sets $\{\alpha_{2j},\alpha_{2j+1}\}$ $(j\geq 0).$
Thus $S_{(B)}$ can be identified with $(\Z/2\Z)^\omega,$
and \mbox{$S_{(B)}\cap S_{(A)}^\mathrm{\,even}$} can be seen to
correspond to the subgroup $\{(a_0,a_1,\ldots)\in(\Z/2\Z)^\omega:
(\forall\,i\geq 0)\;\;a_{2i}=a_{2i+1}\}.$

Let us now choose from the union of the singleton members of $A$
infinitely many elements, which we will denote $\alpha_i$ for $i<0,$
and let $f\in S$ be any permutation such that $\alpha_i f=\alpha_{i+2}$
for all $i\in\Z.$
Then we see that the conjugation map $g\mapsto f^{-1}g\<f$ will carry
$S_{(B)}$ into itself, by a homomorphism which, identifying $S_{(B)}$
with $(\Z/2\Z)^\omega,$ takes the form $(a_0,a_1,\ldots)\mapsto
(0,a_0,a_1,\ldots).$
Now it is not hard to see that every member of $(\Z/2\Z)^\omega$ can be
written (uniquely) as the sum of an element whose $2i\!\<$th and
$2i{+}1\!\<$st coordinates are equal for each $i,$ and an element
whose $0\!\<$th coordinate is $0$ and whose $2i{+}1\!\<$st and
$2i{+}2\!\<$nd coordinates are equal for all $i.$
Hence

\begin{xlist}\item\label{x.SB_fr_even}
$S_{(B)}\;=\;(S_{(B)}\cap S_{(A)}^\mathrm{\,even})\,
(\,f^{-1}(S_{(B)}\cap S_{(A)}^\mathrm{\,even})f\,)\;\leq\;
\langle\<S_{(A)}^\mathrm{\,even}\cup\{f\}\<\rangle.$
\end{xlist}
We also see that $S_{(A)}=S_{(B)}S_{(A)}^\mathrm{\,even}.$
(For, given any $h\in S_{(A)}$ which we wish to represent in
this way, a factor in $S_{(B)}$ can be chosen which
gives a permutation of the desired parity on
each $\Sigma_i\in A,$ and a factor in $S_{(A)}^\mathrm{\,even}$
then turns this into the desired permutation $h.)$
Hence $S_{(A)}\leq\langle\<S_{(A)}^\mathrm{\,even}\cup\{f\}\<\rangle,$
so~(\ref{x.even=all}) holds.

We can now obtain our analog of~(\ref{x.M&N.full}).
Suppose that

\begin{xlist}\item\label{x.likefull}
$H$ is a subgroup of $S,$ and $A$ an infinite family of disjoint
nonempty subsets of $\Omega,$ of unbounded finite cardinalities,
such that writing $\Delta=\bigcup_{\Sigma\in A}\Sigma,$
every member of $\Sym(\Delta)_{(A)}$ extends to an
element of $H_{\{\Delta\}}.$
\end{xlist}
That is, we assume
we can find elements of $H$ which give any specified family of
permutations of the sets $\Sigma$ comprising $A$ -- but we don't
assume that we can control what they do off those sets.
We claim that by adjoining to $H$ one element from $S,$ we
can get a group which contains a subgroup
$S_{(A')}^\mathrm{\,even}$ for some $A'\in\mathcal P$
satisfying the conditions stated before~(\ref{x.even=all}) (infinitely
many singletons, all other members having cardinality $\geq 4).$

To do this let us split the set $A$ of~(\ref{x.likefull})
into three infinite disjoint subsets,
$A=A_1\<\cup\<A_2\<\cup\<A_3,$ in any way such that $A_1$
has members of unbounded finite
cardinalities and no members of cardinality $<4.$
If we let $\Delta_1=\bigcup_{\Sigma\in A_1}\Sigma,$
$\Delta_2=\bigcup_{\Sigma\in A_2}\Sigma,$
$\Delta_3=\Omega-\Delta_1-\Delta_2,$
we see that these sets each have cardinality $\aleph_0$
(the last because it contains $\bigcup_{\Sigma\in A_3}\Sigma).$
Since $\Delta_1\cup\Delta_2\subseteq\Delta,$
it follows from~(\ref{x.likefull}) that

\begin{xlist}\item\label{x.any1_id2}
For every $f\in\Sym(\Delta_1)_{(A_1)}$ there exists an
$f'\in H$ which agrees with $f$ on $\Delta_1,$
and acts as the identity on $\Delta_2.$
\end{xlist}
Now take any $g\in S$ that interchanges $\Delta_2$
and $\Delta_3,$ and fixes $\Delta_1$ pointwise.
Conjugating~(\ref{x.any1_id2}) by $g$ gives

\begin{xlist}\item\label{x.any1_id3}
For every $f\in\Sym(\Delta_1)_{(A_1)}$ there exists an
$f'\in g^{-1}H\,g$ which agrees with $f$ on $\Delta_1$
and acts as the identity on $\Delta_3.$
\end{xlist}

Now if one forms the commutator of a permutation which acts as the
identity on $\Delta_2$ and preserves $\Delta_3$ with
a permutation which acts as the
identity on $\Delta_3$ and preserves $\Delta_2,$
one gets an element which acts as the
identity on $\Delta_2\cup\Delta_3.$
Hence from~(\ref{x.any1_id2}) and~(\ref{x.any1_id3}) we may
conclude that $\langle H\cup\{g\}\<\rangle$ contains elements
which act as the identity on $\Delta_2\cup\Delta_3=\Omega-\Delta_1,$
while acting on each $\Sigma\in A_1$ by any specified commutator
in $\Sym(\Sigma).$
Moreover, in the symmetric group on a finite set $\Sigma,$ the
commutators are precisely the even permutations
\cite[Theorem~1]{OO}; so letting $A'$ be the
partition of $\Omega$ consisting of the members of $A_1$
and all singleton subsets of $\Omega-\Delta_1,$ we have

\begin{xlist}\item\label{x.SA'_even}
$\langle H\cup\{g\}\<\rangle\;\geq\;S_{(A')}^\mathrm{\,even}.$
\end{xlist}
Combining with~(\ref{x.even=all}), we get our analog
of~(\ref{x.M&N.full}), namely

\begin{xlist}\item\label{x.fin.full}
If $H$ satisfies~(\ref{x.likefull}), then
$H\sucq S_{(A')}$ for some (hence by~(\ref{x.SAeqSB}), for all)
$A'\in\mathcal P.$
\end{xlist}

With the help of~(\ref{x.fin.full}) we can now prove a strengthening
thereof, analogous to Lemma~\ref{L.D_i} of the preceding section:

\begin{lemma}\label{L.D_i2}
Let $\Omega$ be a countably infinite set and $G$ a subgroup of
$\Sym(\Omega),$ and suppose there exist a sequence
of distinct elements $(\alpha_i)_{i\in\omega}\in\Omega^\omega,$
an unbounded sequence of positive integers $(N_i)_{i\in\omega},$ and
a sequence of sets $D_i\subseteq\Omega^i$ $(i\in\omega),$ such that
\\[2pt]
{\rm(i)}~ For each $i\in\omega$ and each
$(\beta_0,\ldots,\beta_{i})\in D_{i+1},$ we have
$(\beta_0,\ldots,\beta_{i-1})\in D_i;$
\\[2pt]
{\rm(ii)}~ For each $i\in\omega$ and each
$(\beta_0,\ldots,\beta_{i-1})\in D_i,$ there exist at least
$N_i$ elements $\beta\in\Omega$ such that
$(\beta_0,\ldots,\beta_{i-1},\beta)\in D_{i+1};$ and
\\[2pt]
{\rm(iii)}~ If $(\beta_i)_{i\in\omega}\in\Omega^\omega$ has the property
that $(\beta_0,\ldots,\beta_{i-1})\in D_i$ for each $i\geq 0,$
then there exists $g\in G$ such that $(\beta_i)=(\alpha_i\,g)$
in $\Omega^\omega.$
\vspace{2pt}

Then $G\sucq S_{(A)}$ for some, equivalently, for all $A\in\mathcal P.$
\end{lemma}\begin{proof}
This will be similar to the proof of Lemma~\ref{L.D_i}, but with two
simplifications and one complication.
The simplifications are, first, that we will not need to handle
simultaneously strings of permutations $(\pi_1,\ldots,\pi_r)$ for
all decompositions of $\{0,\ldots,i{-}1\}$ as
$\{0,\ldots,n_1{-}1\}\cup\ldots\cup\{n_{r-1},\ldots,n_r{-}1\},$ but
only for a single decomposition, and, secondly, that we will not have
infinite families of choices that have to be spread out over
successive rounds of the construction, as
discussed at the end of the last section.
The complication is that in general not all of the $N_i$ in our
hypothesis will be large enough for our immediate purposes;
hence each time we move to longer strings of
indices, we will have to jump forward to a value $i=i(j)$ such
that $N_{i(j)}$ is large enough.

We begin by fixing an arbitrary increasing sequence of natural
numbers $0=n_0<n_1<\ldots\<,$ such that
the successive differences $n_m-n_{m-1}$ are unbounded.
We shall now construct recursively integers
$-1=i(-1)<i(0)<\ldots<i(j)<\ldots\<,$ and for
each $r\geq 0$ a subset $E_r\subseteq D_{i(n_r-1)+1}.$
The elements of each $E_r$ will be denoted
$e(\pi_1,\ldots,\pi_r),$ where $\pi_m\in\Sym(\{n_{m-1},n_{m-1}{+}1,
\ldots,n_m{-}1\})$ for $m=1,\ldots,r.$

We again begin with $E_0=D_0=$ the singleton consisting of the
empty string.
Now assume inductively for some $r$ that $i(0),\ldots,i(n_{r-1}{-}1)$
and $E_0,\ldots,E_{r-1}$ have been constructed.
We want to choose $n_r-n_{r-1}$
values $i(n_{r-1}),\ldots,\linebreak[3]i(n_r{-}1)\in\omega$ and 
extend each $e(\pi_1,\ldots,\pi_{r-1})\in E_{r-1}$ to a family
of elements $e(\pi_1,\ldots,\pi_r)\in D_{i(n_r-1)+1},$ obtaining
one such extension for each 
$\pi_r\in\Sym(\{n_{r-1},\ldots,n_r{-}1\}),$ in such a way that

\begin{xlist}\item\label{x.distinct}
The components of each $(i(n_r{-}1)\<{+}\<1)\!\<$-tuple
$e(\pi_1,\,\ldots,\pi_r)$ which correspond to the $n_r-n_{r-1}$ indices
$i(n_{r-1}),\linebreak[0]\;i(n_{r-1}{+}1),\;\ldots,\;i(n_r{-}1)$
are distinct from each other, from those components of the
other $i(n_r)\!\<$-tuples
\linebreak[3] $e(\pi_1',\ldots,\linebreak[0]\pi_r')$
$((\pi_1',\ldots,\linebreak[0]\pi_r')\neq(\pi_1,\ldots,\pi_r))$
corresponding to any of the indices
$i(n_{r-1}),\;i(n_{r-1}+1),\;\ldots,\linebreak[0]\;i(n_r{-}1),$
and also from the components of the elements of~$E_{r-1}$
with indices $i(0),\linebreak[0]\;i(1),\;\ldots,\linebreak[0]
\;i(n_{r-1}{-}1).$
\end{xlist}
Hence let us choose values $i(n_{r-1}),\ldots,i(n_r{-}1)$
such that $N_{i(n_{r-1})},\,N_{i(n_{r-1}+1)},\linebreak[0]
\,\ldots\,,\linebreak[0]
\,N_{i(n_r-1)}$ are all $\geq
|E_{r-1}|\<((n_r\<{-}\<n_{r-1})\,(n_r\<{-}\<n_{r-1})!\,+\,n_{r-1}).$
(The factor $n_r-n_{r-1}$ represents the number of new components
of each string referred to in~(\ref{x.distinct});
$(n_r-n_{r-1})!$ is the number of values of $\pi_r,$ and the final
summand $n_{r-1}$ is the number of components of each
member of $E_{r-1}$ that we also have to avoid.)
Using these $i(j),$ it is not hard to see from our hypothesis~(ii)
that we can indeed extend our strings
$e(\pi_1,\ldots,\pi_{r-1})\in D_{i(n_{r-1}-1)+1}$ to strings
$e(\pi_1,\ldots,\pi_r)\in D_{i(n_r-1)+1}$
so that~(\ref{x.distinct}) holds.

As in the proof of Lemma~\ref{L.D_i} we now choose a
single permutation $s$ of $\Omega,$ this time such that

\begin{xlist}\item\label{x.*b_j.s2}
For each $e(\pi_1,\ldots,\pi_r)=(\beta_j)_{0\leq j<n_r}\in E_r,$
the element $s$ acts on the components $\beta_{i(n_{r-1})},\linebreak[0]
\beta_{i(n_{r-1}+1)},\ldots,\beta_{i(n_r-1)}$ of this tuple
so that $\beta_{i(j)}\<s = \beta_{i(j\pi_r)}.$
\end{xlist}

Now let $\Delta$ be the set $\{\alpha_{i(j)}:j\geq 0\},$
and let $A$ be the partition of $\Delta$ into subsets
$\Sigma_r=\{\alpha_{i(j)}:n_{r-1}\leq j< n_r\}$ $(r\,{\geq}\<1).$
Thus the general element of $\Sym(\Delta)_{(A)}$ has the
form $\alpha_{i(j)}\mapsto\alpha_{i(j\pi)}$ for
some $\pi\in\Sym(\omega)$ that preserves each set
$\{n_{r-1},\ldots,n_r{-}1\}.$
We claim that for any such permutation $\pi,$ there is a $g\in G$ such
that $s$ ``acts as $\pi$ on the subscripts'' of the translated sequence
$(\alpha_{i(j)}\,g)_{\strt j\geq 0},$ i.e., such that for all $j\geq 0,$

\begin{xlist}\item\label{x.*a_i.g.s2}
$\alpha_{i(j)}\,g\,s\;=\;\alpha_{i(j\pi)}\,g.$
\end{xlist}

Indeed, given $\pi,$ if for each $r\geq 1$ we let
$\pi_r\in\Sym(\{n_{r-1},\ldots,n_r{-}1\})$ denote the restriction
of $\pi$ to $\{n_{r-1},\ldots,\linebreak[0]n_r{-}1\},$
then as in the proof of Lemma~\ref{L.D_i}, the strings
$e(\<),$ $e(\pi_1),$ $e(\pi_1,\pi_2),\ldots$ fit together to
give a string $(\beta_i)$ such that~(\ref{x.*b_j.s2}) says that
$s$ ``acts like $\pi$'' on the components of $(\beta_i)$ indexed by the
$i(j)$ $(j\in\omega).$
By hypothesis~(iii), we can write $(\beta_i)$ as $(\alpha_i\,g)$
so this condition becomes~(\ref{x.*a_i.g.s2}).

But~(\ref{x.*a_i.g.s2}) can be read as saying that $g\,s\,g^{-1}$ acts
on $\Delta$ as the arbitrary element
$\alpha_{i(j)}\mapsto\alpha_{i(j\pi)}$ of $\Sym(\Delta)_{(A)};$
hence letting $H=\langle\<G\cup\{s\}\rangle,$ (\ref{x.likefull})
holds, so by~(\ref{x.fin.full}), $G\sucq S_{(A')}$
for some $A'\in\mathcal P,$ completing the proof of the lemma.
\end{proof}

The next argument also parallels what we did
in the preceding section (though it will be less convoluted):
For any $G\leq\Sym(\Omega)$ satisfying~(\ref{x.unbdd_orbs}), we shall
obtain families $D_i$ satisfying conditions~(i) and~(ii) of the above
lemma, and such that if $G$ is closed, condition~(iii) also holds.

Assume $G\leq\Sym(\Omega)$ satisfies~(\ref{x.unbdd_orbs}), and fix an
unbounded sequence of positive integers $(N_i)_{i\in\omega}.$
We shall begin by constructing for each $j\geq 0$ a certain element
$\alpha_j,$ and a certain finite subset $K_j$ of $G,$
which will be indexed

\begin{xlist}\item\label{x.g()2}
$K_j\;=\;\{\<g(k_0,k_1,\ldots,k_{j-1})\;:\;
0\leq k_i<N_i\;\;(0\leq i<j)\<\}.$
\end{xlist}
Again, $K_0$ will have only one member,
$g(\<),$ which we take to be $1\in G.$

Let us assume inductively for some $j\geq 0$ that elements
$\alpha_i$ have been defined for all $i<j$ and that subsets
$K_i$ have been defined for all $i\leq j.$
Let $\Gamma_j\subseteq\Omega$ denote (essentially as before)
the set of images of $\eps_0,\ldots,\eps_{j-1}$ and of
$\alpha_0,\ldots,\alpha_{j-1}$ under inverses
of elements of $K_0\cup\ldots\cup K_j.$
Let $\alpha_j$ be any element of $\Omega$ not fixed by
$G_{(\Gamma_j)},$ whose orbit under that group has cardinality
at least $|K_j|\,N_j;$ such an element exists by (\ref{x.unbdd_orbs}).

We now fix an arbitrary order in which we shall construct
the elements $g(k_0,k_1,\linebreak[3]\ldots,k_j)$ of $K_{j+1}.$
When it is time to construct $g(k_0,k_1,\ldots,k_j),$
we set $g'=g(k_0,k_1,\ldots,k_{j-1}),$ and
left-multiply this by any element $h\in G_{(\Gamma_j)}$
with the property that $\alpha_j\,h\<g'$ is distinct from the
images of $\alpha_j$ under those elements of $K_j$ so far constructed.
Our choice of $\alpha_j$ insures that its orbit under $G_{(\Gamma_j)}$
is large enough so that collisions with all such elements can be
avoided, and we define $g(k_0,k_1,\ldots,k_j)$ to be the
product $h\<g'.$

For each $i$ we then define the sets $D_i$ by

\begin{xlist}\item\label{x.def_D_i2}
$D_i\;=\;\{\<(\alpha_0\<g,\<\ldots\<,\<\alpha_{i-1}\<g)\,:\,
g\in K_i\<\}.$
\end{xlist}

We now see exactly as before that conditions~(i)
and~(ii) of Lemma~\ref{L.D_i2} are satisfied, and that if $G$ is
closed, we can use Lemma~\ref{L.b+f} to get condition~(iii) as well,
so by Lemma~\ref{L.D_i2}, $G\sucq S_{(A)}$ for some $A\in\mathcal P.$

On the other hand, the reverse inequality is immediate:
Taking any $\Gamma$ as in the first clause of~(\ref{x.unbdd_orbs})
and letting $B$ denote the set of orbits of $G_{(\Gamma)},$ so that
$B\in\mathcal P,$
we get $G\approx G_{(\Gamma)}\leq S_{(B)}\approx S_{(A)}$
(where the first relation holds by Lemma~\ref{L.G_*G} and
Lemma~\ref{L._0<=>_1}(ii)).
Combining these inequalities we have $G\approx S_{(A)}.$

This completes the main work of the proof of

\begin{theorem}\label{T.unbdd_orbs}
Let $\Omega$ be a countably infinite set,
and $\mathcal P$ the set of partitions of $\Omega$
defined in~{\rm(\ref{x.def_P})}.
Then the subgroups $S_{(A)}\leq S$ with $A\in\mathcal P$ {\rm(}which
are clearly all closed\/{\rm)} are mutually $\approx\!\<$-equivalent,
and a closed subgroup $G\leq S$ belongs to the equivalence class of
those subgroups if and only if it satisfies~{\rm(\ref{x.unbdd_orbs})}.

Moreover, the members of this $\approx\!\<$-equivalence class are
$\prec$ the members of the
equivalence class of Theorem~\ref{T.inf_orbs}.
\end{theorem}\begin{proof}

We have so far proved mutual equivalence of the $S_{(A)},$
and the sufficiency of~(\ref{x.unbdd_orbs}) for membership of a closed
subgroup $G$ in their common equivalence class.
To see necessity, consider any closed subgroup
$G$ which does not satisfy~(\ref{x.unbdd_orbs}).
Then either $G$ satisfies~(\ref{x.inf_orbs}),
or there exists a finite set $\Gamma$
such that $G_{(\Gamma)}$ has orbits of bounded finite cardinality.

In the former case, Theorem~\ref{T.inf_orbs} shows that
$G\approx S;$ but from the ``only if'' direction of that theorem
we see that for $A\in\mathcal P$ we have $S_{(A)}\not\approx S,$ and
hence $G\not\approx S_{(A)}.$

In the case where some
$G_{(\Gamma)}$ has all orbits of bounded finite cardinality,
let $A$ be the partition of $\Omega$ consisting of those orbits.
Then $G\approx G_{(\Gamma)}\preq S_{(A)},$ and by the
last sentence of Theorem~\ref{T.no_uncrwdd}, $S_{(A)}$ is not $\sucq$
the members of the equivalence class of this section,
hence $G$ is not in that equivalence class.

In the final sentence, the inequality $\preq$ holds because the
equivalence class of Theorem~\ref{T.inf_orbs} contains $S$ itself.
We have just seen that the two classes in question are distinct, so we
have strict inequality $\prec.$
\end{proof}

\section{Orbits of bounded size.}\label{S.bdd_orbs}
Moving on to still smaller subgroups,
we now consider $G\leq S$ satisfying

\begin{xlist}\item\label{x.bdd_orbs}
There exists a finite subset $\Gamma\subseteq\Omega$ and a positive
integer $n$ such that the cardinalities of all the orbits of
$G_{(\Gamma)}$ are bounded by $n,$ but there exists no
such $\Gamma$ with $G_{(\Gamma)}=\{1\}.$
\end{xlist}

Analogously to~(\ref{x.def_P}), we define

\begin{xlist}\item\label{x.def_Q}
$\mathcal Q\;=\;\{A: A$ is a partition of $\Omega$ for which there is
a common finite bound to the cardinalities of the members of $A,$ and
such that infinitely many members of $A$ have cardinality $>1\}.$
\end{xlist}

Unlike the $\mathcal P$ of the preceding section, $\mathcal Q$
has, up to isomorphism, a natural distinguished member, namely
a least isomorphism class with respect to refinement:

\begin{xlist}\item\label{x.A0}
We will denote by $A_0$ an element of $\mathcal Q,$ unique up
to isomorphism, which has infinitely many $1\!\<$-element
members, infinitely many $2\!\<$-element members, and no others.
\end{xlist}

Clearly any $A\in\mathcal Q$ can be refined to a partition $A'_0$
isomorphic to $A_0,$ hence $S_{(A)}\geq S_{(A'_0)}\approx S_{(A_0)},$
so $S_{(A)}\sucq S_{(A_0)}.$
We claim that the reverse inequality
$S_{(A)}\preq S_{(A_0)}$ also holds.
To show this, let us draw a graph with the elements
of $\Omega$ as vertices, and with edges making each member of our
given partition $A$ a chain (in an arbitrary way), and no other edges.
Now color the edges of each such chain alternately red and green,
subject to the condition that infinitely many chains have a
terminal red edge and infinitely many have a terminal green edge.
Clearly, the partition of $\Omega$ whose non-singleton members
are the pairs of points linked by red edges, all other
points forming singletons, is isomorphic to $A_0;$
hence the group of permutations whose general member acts by
transposing an arbitrary subset of the red-linked pairs of vertices
and fixing everything else can
be written $f^{-1}S_{(A_0)}f$ for some $f\in\Sym(\Omega).$
Similarly, the group of permutations which act by
transposing some pairs of green-linked vertices
and fixing everything else can be written $g^{-1}S_{(A_0)}\<g.$
Moreover, for each $\Sigma\in A,$ any permutation of $\Sigma$ can
be obtained by composing finitely many permutations, each of which
acts either by interchanging only red-linked pairs or
by interchanging only green-linked pairs (this is easiest to see
by looking at permutations that interchange one such pair at a time);
and the number of such
factors needed can be bounded in terms of the cardinality of $\Sigma\<.$
Since there is a common bound to the cardinalities of
the sets $\Sigma\in A,$ we see that every member of
$S_{(A)}$ can be written as a finite
product of members of $f^{-1}S_{(A_0)}f$ and $g^{-1}S_{(A_0)}g,$ so
$\langle\<S_{(A_0)}\cup\{f,g\}\<\rangle\geq S_{(A)},$
so $S_{(A_0)}\sucq S_{(A)}.$
Combining this with the observation
at the start of this paragraph, we get $S_{(A_0)}\approx S_{(A)},$ so

\begin{xlist}\item\label{x.SAeqSB2}
$S_{(A)}\;\approx\;S_{(B)}$~ for all $A,B\in\mathcal Q\<.$
\end{xlist}

We obtain next the result that will play the role
that~(\ref{x.M&N.full}) played in~\S\ref{S.inf_orbs}
and~(\ref{x.fin.full}) played in~\S\ref{S.unbdd_orbs}.
The development will be similar to the latter case, though simpler.
Suppose that
\begin{xlist}\item\label{x.likefull2}
$H$ is a subgroup of $S,$ and $A$ an infinite family of disjoint
$2\!\<$-element subsets of $\Omega$ such that,
writing $\Delta=\bigcup_{\Sigma\in A}\Sigma,$ every member
of $\Sym(\Delta)_{(A)}$ extends to an element of $H_{\{\Delta\}}.$
\end{xlist}
(Again we do not assume we have any control over the behavior of these
elements outside of $\Delta,$ though again our goal will be to
get such control in an extended subgroup.)
Let us index $A$ by $\Z,$ writing $A=\{\Sigma_i:i\in\Z\},$
and let $h$ be an element of $S$ which for each $i\in\Z$
sends $\Sigma_i$ bijectively to $\Sigma_{i+1},$ and
which fixes all elements of $\Omega-\Delta\<.$
We claim that as $f$ runs over all elements of $H_{\{\Delta\}}$
that extend elements of $\Sym(\Delta)_{(A)},$
the commutators $h^{-1}f^{-1}h\<f$ all fix $\Omega-\Delta$
pointwise, and their restrictions to $\Delta$ give all
elements of $\Sym(\Delta)_{(A)}.$
The first fact holds because $h$ fixes $\Omega-\Delta$ pointwise.
The second may be seen by looking at $h^{-1}f^{-1}h\<f$ as
$(h^{-1}f^{-1}h)\<f,$
noting that both factors are members of $\Sym(\Delta)_{(A)},$
and examining their behaviors on the general $2\!\<$-element set
$\Sigma_i\in A.$
One sees that $(h^{-1}f^{-1}h)\<f$ acts by
the trivial permutation on $\Sigma_i$ if and only if $f$ acts trivially
either on both of $\Sigma_{i-1}$ and $\Sigma_i,$ or on neither,
while $(h^{-1}f^{-1}h)\<f$ acts by the nonidentity element of
$\Sym(\Sigma_i)$ in the remaining cases.
One easily deduces that by appropriate choice of $f$ one can get
an arbitrary action on the family of subsets $\Sigma_i.$

Thus $\langle H\cup\{h\}\rangle$ contains
a subgroup conjugate in $S$ to $S_{(A_0)},$ proving

\begin{xlist}\item\label{x.2-full}
If $H$ satisfies~(\ref{x.likefull2}), then $H\sucq S_{(A_0)}.$
\end{xlist}

The result analogous to Lemmas~\ref{L.D_i}
and~\ref{L.D_i2} will be quite simple to state and prove this time:

\begin{lemma}\label{L.D_i3}
Let $\Omega$ be a countably infinite set and $G$ a subgroup of
$\Sym(\Omega),$ and suppose there exist two disjoint sequences
of distinct elements, $(\alpha_i),\,(\beta_i)\in\Omega^\omega,$
such that for every element $(\gamma_i)\in
\prod_{i\in\omega}\{\alpha_i,\beta_i\}\subseteq\Omega^\omega,$
there exists $g\in G$ such that $(\gamma_i)=(\alpha_i\,g).$

Then $G\sucq S_{(A_0)}.$
\end{lemma}\begin{proof}

Let $\Delta=\{\alpha_i:i\in\omega\},$
let $A$ be the partition of $\Delta$ whose members are the two-element
sets $\{\alpha_{2j},\<\alpha_{2j+1}\}$ $(j\geq 0),$ and let
$s\in\Sym(\Omega)$ be any element which fixes all the elements
$\alpha_i$ and interchanges $\beta_{2j}$ and $\beta_{2j+1}$
for all $j\geq 0.$
We claim that every member of $\Sym(\Delta)_{(A)}$
extends to an element of $\langle\<G\cup\{s\}\rangle.$

Indeed, given $f\in\Sym(\Delta)_{(A)},$ define
$(\gamma_i)\in\Omega^\omega$
by letting $\gamma_i=\beta_i$ if $\alpha_i$ is moved
by $f$ (i.e., if it is transposed with the other member
of its $A\!\<$-equivalence class) and $\gamma_i=\alpha_i$ otherwise.
By hypothesis we can find $g\in G$ such that $\gamma_i=\alpha_i\,g$ for
all $i.$
It is now easy to see that $g\,s\,g^{-1}$ acts by $f$ on $\Delta\<.$

Thus $\langle\<G\cup\{s\}\rangle$ satisfies~(\ref{x.likefull2}), so
by~(\ref{x.2-full}),
$S_{(A_0)}\preq\langle\<G\cup\{s\}\rangle\approx G.$
\end{proof}

As the pattern of the two preceding sections suggests, we will
now prove that any closed subgroup $G$ satisfying~(\ref{x.bdd_orbs})
satisfies the hypothesis of the above lemma.
We begin with a reduction: Assuming~(\ref{x.bdd_orbs}),
let $M>1$ be the largest integer such that for every finite
subset $\Gamma\subseteq\Omega,$ the group $G_{(\Gamma)}$ has
orbits of cardinality at least $M.$
Thus, there exists some finite $\Delta$ such that $G_{(\Delta)}$
has no orbits of cardinality $>M.$
Since $G_{(\Delta)}$
inherits from $G$ the property~(\ref{x.bdd_orbs}), we may replace $G$
by $G_{(\Delta)}$ and so assume without loss of generality that

\begin{xlist}\item\label{x.M=max}
For every finite subset $\Gamma\subseteq\Omega,$ the maximum
of the cardinalities of the orbits of $G_{(\Gamma)}$ is $M.$
\end{xlist}
A consequence is that for any such $\Gamma,$ every orbit
of $G_{(\Gamma)}$ of cardinality $M$ is also an orbit of $G$
(since the orbit of $G$ containing it cannot have larger cardinality).
Thus

\begin{xlist}\item\label{x.orbtimesG}
If $\Gamma$ is a finite subset of $\Omega,$ and $\alpha$ an
element of $\Omega$ such that $|\<\alpha\,G_{(\Gamma)}\<|=M,$
then for every $g\in G$ we have
$\alpha\,G_{(\Gamma)}\,g=\alpha\,G_{(\Gamma)}.$
\end{xlist}

We shall now construct recursively, for each $j\geq 0,$
elements $\alpha_j,\,\beta_j\in\Omega$ and a subset
$K_j\subseteq G,$ indexed as

\begin{xlist}\item\label{x.g()3}
$K_j\;=\;\{\<g(k_0,k_1,\ldots,k_{j-1})\;:\;
(k_0,k_1,\ldots,k_{j-1})\in\{0,1\}^j\<\}.$
\end{xlist}

Again we start with $K_0=\{g(\<)\}=\{1\}.$
Assuming inductively for some $j\geq 0$ that $\alpha_i,\ \beta_i$
have been defined for all $i<j$ and $K_i$ for all $i\leq j,$
we let $\Gamma_j\subseteq\Omega$ denote
the set all of images of $\eps_0,\ldots,\eps_{j-1},$
$\alpha_0,\ldots,\alpha_{j-1},$
$\beta_0,\ldots,\beta_{j-1}$ under inverses
of elements of $K_0\cup\ldots\cup K_j.$
By assumption, $G_{(\Gamma_j)}$ has an $M\!$-element orbit.
Let $\alpha_j$ and $\beta_j$ be any
two distinct elements of such an orbit.
(Note that $\alpha_j$ and $\beta_j$ are distinct from all
$\alpha_i,$ $\beta_i$ for $i<j,$
since the latter are fixed by $G_{(\Gamma_j)}.)$
For each $(k_0,\ldots,k_{j-1})\in\{0,1\}^j,$ we let
$g(k_0,\ldots,k_{j-1},0)$ and
$g(k_0,\ldots,k_{j-1},1)$ be elements of $G$ obtained
by left-multiplying $g(k_0,\ldots,k_{j-1})\in K_j$
by an element $h\in G_{(\Gamma_j)},$ chosen so that
$\alpha_j\<h\<g(k_0,\ldots,k_{j-1})$ is
$\alpha_j,$ respectively $\beta_j.$
This is possible by~(\ref{x.orbtimesG}).

Given an infinite string $(k_i)$ of $0\!\<$'s and $1\!\<$'s,
the elements $g(k_0,k_1,\ldots,k_{j-1})$
will again converge in $S$ by Lemma~\ref{L.b+f}.
Assuming $G$ closed, the limit belongs to $G,$ and
clearly gives us the hypothesis of Lemma~\ref{L.D_i3},
hence the conclusion that $G\sucq S_{(A_0)}.$

Again we easily get the reverse inequality: Taking
$\Gamma$ as in the first clause of~(\ref{x.bdd_orbs})
and letting $A$ denote the partition of $\omega$ into orbits
of $G_{(\Gamma)},$ we have $G\approx G_{(\Gamma)}\leq
S_{(A)}\approx S_{(A_0)}$ by~(\ref{x.SAeqSB2}).

We deduce

\begin{theorem}\label{T.bdd_orbs}
Let $\Omega$ be a countably infinite set,
and $\mathcal Q$ the set of partitions of $\Omega$
defined in~{\rm(\ref{x.def_Q})}.
Then the subgroups $S_{(A)}\leq S$ for $A\in\mathcal Q$ {\rm(}which
are clearly closed\/{\rm)} are mutually $\approx\!\<$-equivalent,
and a closed subgroup $G\leq S$ belongs to the equivalence class of
those subgroups if and only if it satisfies~{\rm(\ref{x.bdd_orbs})}.

The members of this $\approx\!\<$-equivalence class are $\prec$ the
members of the equivalence class of Theorem~\ref{T.unbdd_orbs}.
\end{theorem}\begin{proof}
This is obtained using the above results exactly
as Theorem~\ref{T.unbdd_orbs} was obtained from the results
of the preceding section, except that we need a different argument
to show that $G$ does not belong to the
$\approx\!\<$-equivalence class in question
if it does not satisfy the final clause of~(\ref{x.bdd_orbs}),
i.e., if there exists a finite subset $\Gamma\subseteq\Omega$
such that $G_{(\Gamma)}=\{1\}.$
In that situation, any subgroup $\approx G$ will be
$\approx\{1\},$ hence countable;
but clearly $S_{(A_0)}$ is uncountable, so $S_{(A_0)}\not\approx G.$
\end{proof}

\section{Countable subgroups.}\label{S.countable}

The final step of our classification is now easy, and we even get a
little extra information:

\begin{theorem}\label{T.countable}
The countable subgroups of $S=\Sym(\Omega)$ form an equivalence
class under ${\approx}\<,$ and members of this class are $\prec$ the
members of the equivalence class of Theorem~\ref{T.bdd_orbs}.
Moreover, for $G\leq S,$ the following conditions are equivalent.
\\[2pt]
{\rm(i)}~ $G$ is countable and closed.
\\[2pt]
{\rm(ii)}~ There exists a finite subset
$\Gamma\subseteq\Omega$ such that $G_{(\Gamma)}=\{1\}.$
\\[2pt]
{\rm(iii)}~ $G$ is discrete.

\end{theorem}\begin{proof}
The countable subgroups are clearly the subgroups
$\approx_{\strt\aleph_1}\!\!\{1\},$
and as noted in Lemma~\ref{L._0<=>_1},
for subgroups of symmetric groups $\Sym(\Omega),$
$\approx_{\strt\aleph_1}\!\!$-equivalence is the same as
$\approx_{\strt\aleph_0}\!\!$-equivalence, which is what we are
calling $\approx\!\<$-equivalence.
This gives the first assertion; the second is also immediate,
since the trivial subgroup is $\preq$ all subgroups,
and is $\not\approx$ the subgroups of Theorem~\ref{T.bdd_orbs} by
the ``only if'' assertion of that theorem.

To prove the equivalence of~(i)--(iii), we note first that~(ii)
and~(iii) are equivalent, since a neighborhood basis of the
identity in the function topology on $G$ is given by the subgroups
$G_{(\Gamma)}$ for finite $\Gamma,$ so the identity element (and hence
by translation, every element) is isolated in $G$ if and only if
some such subgroup is trivial.

To see that these equivalent conditions imply~(i), observe
that~(ii) implies that $G\approx G_{(\Gamma)}=\{1\},$ hence that
$G$ is countable, while~(iii) implies that $G$ is closed, by general
properties of topological groups.
(If $G$ is a discrete subgroup of a topological group $S,$
take a neighborhood $U$ of $1$ in $S$ containing no
nonidentity element of $G,$ and then a neighborhood $V$
of $1$ such that $VV^{-1}\subseteq U.$
One finds that for any $x\in S,$ $xV$ is a
neighborhood of $x$ containing at most one element of $G;$
so $G$ has no limit points in $S.)$

Conversely, we have seen that any countable $G$ is $\prec$ the members
of the equivalence class of Theorem~\ref{T.bdd_orbs}, hence does not
belong to the equivalence class of any of Theorems~\ref{T.inf_orbs},
\ref{T.unbdd_orbs} or~\ref{T.bdd_orbs}.
Hence if $G$ is also closed, those theorems exclude all
possible behaviors of its subgroups $G_{(\Gamma)}$ (for $\Gamma$ finite)
other than that there exist such a $\Gamma$ with $G_{(\Gamma)}=\{1\};$
so~(i) implies~(ii).
\end{proof}

For convenience in subsequent discussion, let us name the four
equivalence classes of subgroups of $S=\Sym(\Omega)$ which we have
shown to contain all closed subgroups:

\begin{xlist}\item\label{x.C1-C0}
$\!\mathcal C_S\;=$ the $\approx\!\<$-equivalence class of $S.$
\\[2pt]
$\mathcal{C_P}\;=$ the $\approx\!\<$-equivalence class to which
$S_{(A)}$ belongs for all $A\in\mathcal P.$
\\[2pt]
$\mathcal{C_Q}\;=$ the $\approx\!\<$-equivalence class to which
$S_{(A)}$ belongs for all $A\in\mathcal Q\<.$
\\[2pt]
$\mathcal C_1\;\;=$ the $\approx\!\<$-equivalence class
consisting of the countable subgroups\\
\indent of $S.$
\end{xlist}

\section{Notes and questions
on groups of bounded permutations.}\label{S.q.FN}

It would be of interest to investigate the equivalence relation
$\approx$ on classes of subgroups $G\leq\Sym(\Omega)$
other than the class of closed subgroups.
One such class is implicit in the techniques used above:
If $\Omega$ is any set and $d$ a generalized metric on $\Omega,$
let us define the subgroup

\begin{xlist}\item\label{x.defFN(O,d)}
$\mathrm{FN}(\Omega,d)\;\;=\;\;\{\<g\in\Sym(\Omega):{||g||}_d<\infty\}.$
\end{xlist}
(We write $\mathrm{FN},$ for ``finite norm'', rather than $B$ for
``bounded'' to avoid confusion with the symbol for an open ball.)
These subgroups are not in general closed.
For instance, if $d$ does not assume the value $\infty$
(i.e., if it is an ordinary metric) but is unbounded
(say the standard distance metric on
$\Omega=\omega\subseteq\mathbb R),$ then by the former condition,
$\mathrm{FN}(\Omega,d)$ contains all permutations of $\Omega$ that
move only finitely many elements, which form a
dense subgroup of $S=\Sym(\Omega),$ while by the unboundedness of $d,$
it is nevertheless a proper subgroup of $S,$ hence it is not closed.

Here are some easy results about the relation $\preq$ on these
subgroups.
(Cf.\ also~\cite{Such90}.)
Below, ``uncrowded'' means $\aleph_0\!$-uncrowded.

\begin{lemma}\label{L.FN}
Let $d$ be a generalized metric on a countably infinite set $\Omega\<.$
\\[2pt]
{\rm(i)}~ If $d$ is not uncrowded,
then $\mathrm{FN}(\Omega,d)\in\mathcal C_S.$
\\[2pt]
{\rm(ii)}~ If $d$ is uncrowded but not uniformly uncrowded, then
$\mathrm{FN}(\Omega,d)$ is $\sucq$ the groups in
$\mathcal{C_P},$ but is $\not\in\mathcal C_S.$
\\[2pt]
{\rm(iii)}~ If $d$ is uniformly uncrowded, but for some $r<\infty,$
infinitely many of the balls $B_d(\alpha,r)$ contain more than one
element, then $\mathrm{FN}(\Omega,d)$ is $\sucq$ the groups in
$\mathcal{C_Q},$ but $\not\sucq$ the groups in $\mathcal{C_P}.$
\\[2pt]
{\rm(iv)}~ If $d$ is uncrowded and for each $r<\infty$ all but finitely
many balls $B_d(\alpha,r)$ are
singletons, then $\mathrm{FN}(\Omega,d)\in\mathcal C_1.$
\end{lemma}\begin{proof}
In situation~(i), let $B_d(\alpha,r)$ be a ball of finite
radius containing infinitely many elements.
Then all $g\in\Sym(\Omega)_{(\Omega-B_d(\alpha,r))}$
satisfy ${||g||}_d\leq 2r,$ hence lie in $\mathrm{FN}(\Omega,d),$
and the conclusion follows by~(\ref{x.M&N.full}).

In cases~(ii) and~(iii) we can similarly find subgroups
of $\mathrm{FN}(\Omega,d)$ of the form $S_{(A)}$ for $A\in\mathcal P,$
respectively $A\in\mathcal Q,$ while the last sentence of
Theorem~\ref{T.bddG} gives the negative statements for these two cases.

For $d$ as in~(iv), each set $\{g\in\Sym(\Omega):{||g||}_d<n\}$
$(n\in\omega)$ is finite, so their
union, $\mathrm{FN}(\Omega,d),$ is countable.
\end{proof}

Thus, if a group $\mathrm{FN}(\Omega,d)$ belongs to one of the
four $\approx\!\<$-equivalence classes of~(\ref{x.C1-C0}),
the above lemma determines precisely which class that must be.

Note that our definition (\ref{x.defFN(O,d)}) can be rewritten

\begin{xlist}\item\label{x.B=cupB_d}
$\mathrm{FN}(\Omega,d)\;=\;\bigcup_{n\in\omega}\,\{\<g\in
S:{||g||}_d<n\}.$
\end{xlist}
We claim that if the generalized metric $d$ is uncrowded, then for each
$n$ the set $\{g\<{\in}\< S:{||g||}_d<n\}$ is compact.
Indeed, the condition $||g||_d<n$ determines,
for each $\alpha\in\Omega,$ a certain finite set of possibilities for
$\alpha\<g\<;$ so $\{\<g\<{\in}\< S:{||g||}_d<n\}$ is the intersection
of $S$ with a certain compact subset of $\Omega^\Omega\<.$
But for each $\alpha\in\Omega,$ the condition $||g||_d<n$ also
limits us to finitely many possibilities for $\alpha\<g^{-1},$
from which it can be deduced that any limit in $\Omega^\Omega$ of
elements $g\in S$ with $||g||_d<n$ is
again surjective, hence again belongs to $S.$
So $\{g\<{\in}\< S:{||g||}_d<n\}$ is closed in a compact subset
of $\Omega^\Omega,$ hence, as claimed, is compact.
This makes $\mathrm{FN}(\Omega,d)$ a countable union of compact sets,
suggesting the second part of

\begin{question}\label{Q.still4}
If $d$ is an uncrowded generalized metric on a countably infinite
set $\Omega,$ must $\mathrm{FN}(\Omega,d)$ belong to one of the
$\approx\!\<$-equivalence classes of~{\rm(\ref{x.C1-C0})}?

More generally, does every subgroup of $\Sym(\Omega)$ that is a
union of countably many compact subsets belong one of these classes?
What about subgroups that are unions of countably many closed subsets?
What about Borel subgroups?
Analytic subgroups?

If the answer to any of these questions is negative, can one describe
all the \mbox{$\approx\!\<$-equivalence} classes to which such
subgroups belong?
\end{question}

Let us sketch a couple of cases where it is not hard to show that
$\mathrm{FN}(\Omega,d)$ does belong to one of the
equivalence classes of~{\rm(\ref{x.C1-C0})}.

Let $d$ be the standard distance function on
$\omega$ (inherited from $\mathbb R).$
To show that $\mathrm{FN}(\omega,d)\in\mathcal{C_Q},$ let $A_1$
be the partition of $\omega$ consisting of the subsets $\{2m,\<2m{+}1\}$
$(m\in\omega),$ and $A_2$ the partition consisting of the subsets
$\{2m{+}1,2m{+}2\}$ and the singleton $\{0\}.$
Then $A_1,A_2\in\mathcal Q,$ so
$\langle\<S_{(A_1)}\cup S_{(A_2)}\rangle\in\mathcal{C_Q}.$
This subgroup is clearly contained in $\mathrm{FN}(\omega,d);$
we claim that equality holds.

Indeed, given $f\in\mathrm{FN}(\omega,d)$ with ${||f||}_d=n,$ if we let
$\Sigma_i=\{ni,ni{+}1,\ldots,\linebreak[3]n(i{+}1)-1\}$ for $i\geq 0$
and $\Sigma_{-1}=\varnothing,$ then we see that for all $i\geq 0,$
$\Sigma_i f\subseteq \Sigma_{i-1}\cup\Sigma_i\cup\Sigma_{i+1}.$
Letting $B_1$ be the partition of $\omega$ into the subsets
$\Sigma_{2i}\cup\Sigma_{2i+1}$ and $B_2$ the partition into the subsets
$\Sigma_{2i-1}\cup\Sigma_{2i}$ $(i\geq 0),$ we see as in the second
paragraph of the proof of Lemma~\ref{L.loc_bdd} that
$f\in S_{(B_1)} S_{(B_2)}.$
On the other hand, it is easy to show that $S_{(B_1)}$ and $S_{(B_2)}$
are both contained in $\langle\<S_{(A_1)}\cup S_{(A_2)}\rangle,$ using
the fact that any permutation of a $2n\!\<$-element string of integers
$\Sigma_{2i}\cup\Sigma_{2i+1}$ or $\Sigma_{2i-1}\cup\Sigma_{2i}$ can be
written as a product of finitely many transpositions of consecutive
terms, and that the number of transpositions needed can be bounded
in terms of $n$ (cf.\ end of paragraph preceding~(\ref{x.SAeqSB2})).
So $f\in\langle\<S_{(A_1)}\cup S_{(A_2)}\rangle,$ as claimed,
so $\mathrm{FN}(\omega,d)\in\mathcal{C_Q}.$

In the above example, the argument cited from the proof of
Lemma~\ref{L.loc_bdd} uses the fact that for any
$f\in \mathrm{FN}(\omega,d),$ the number of elements that
$f$ carries upward past a given point is equal to the
number that it carries downward past that point.
If we modify this example by replacing $\omega$ with $\Z,$ again with
the standard metric, that property no longer holds, as shown by the
translation function $t\,{:}\;n\mapsto n{+}1.$
It is not hard to see, however, that given $f\in \mathrm{FN}(\Z,d),$
the difference between the number of elements that $f$ moves upward
and downward past a given point is the same for all points, and that
the function associating to $f$ the common value of this difference
is a homomorphism $v\<{:}\;\mathrm{FN}(\Z,d)\to\Z.$
If we let $A_1$ denote the partition of $\Z$ into sets $\{2m,\<2m{+}1\}$
and $A_2$ the partition into sets $\{2m{+}1,2m{+}2\},$ we see that
$S_{(A_1)}$ and $S_{(A_2)}$ lie in the kernel of $v,$ while $v(t)=1.$
The argument of the preceding paragraph can be adapted to show
that $\langle\<S_{(A_1)}\cup S_{(A_2)}\<\rangle=\ker(v),$ hence that
$\langle\<S_{(A_1)}\cup S_{(A_2)}\cup\{t\}\<\rangle=\mathrm{FN}(\Z,d);$
so this group also belongs to~$\mathcal{C_Q}.$
(For some further properties of this example see
Suchkov \cite{Such84}, \cite{Such85},
where $\mathrm{FN}(\Z,d)$ and its subgroup
$\mathrm{ker}(t)$ are called $\bar G$ and $G$ respectively.)

An example that falls under case~(ii) of Lemma~\ref{L.FN} (so
that if $\mathrm{FN}(\Omega,d)$ belongs to one of our four classes,
that class is $\mathcal{C_P})$ is given by
$\Omega=\{\sqrt n:n\in\omega\}$ with the metric induced from
$\mathbb R\<.$
We suspect one can show that it does belong to $\mathcal{C_P}$
by adapting the method we used for $\mathrm{FN}(\omega,d),$
putting in the roles of $A_1$
and $A_2$ the partitions of $\Omega$ arising from the integer-valued
functions $\alpha\mapsto[\alpha/2]$ and
$\alpha\mapsto[(\alpha{+}1)/2],$
where $[-]$ denotes the integer-part function.

Two cases that have some similarity to that of
$\mathrm{FN}(\Z,d)$ but seem less trivial, and might
be worth examining, are those given by the vertex-sets
of the Cayley graphs of the free abelian group, respectively the free
group, on two generators, with the path-length metric.
An example of a different sort is the set $\omega$ with the ultrametric
under which $d(\alpha,\beta)$ is the greatest $n$ such that $\alpha$
and $\beta$ differ in the $\!\<n\!\<$th digit of their
base-$\!2$ expansions.
 From the fact that this $d$ is an uncrowded {\em ultra}\<metric,
it is easily deduced that $\mathrm{FN}(\omega,d)$ is
the union of a countable chain of compact sub{\em groups}.
All three of these examples fall under case~(iii) of Lemma~\ref{L.FN},
so that if they belong to any of the classes of~(\ref{x.C1-C0})
it is $\mathcal{C_Q}.$

\section{Further questions about $\preq$
and ${\approx}.$}\label{S.q.preq}
It seems unlikely that one can in any reasonable sense describe
{\em all} $\approx\!\<$-equivalence classes of subgroups of the
symmetric group on a countably infinite set $\Omega\<.$
On the other hand, if one regards the set of such equivalence classes
as a join-semilattice, with join operation induced by the
map $(G_1,\,G_2)\mapsto\langle\<G_1\cup G_2\<\rangle$ on subgroups,
one may ask about the properties of this semilattice.
The cardinal $|\<G\<|+\aleph_0$ is an $\approx\!\<$-invariant on
subgroups $G$ of $\Sym(\Omega),$ and induces a homomorphism from this
join-semilattice onto the semilattice of cardinals between $\aleph_0$
and $2^{\aleph_0}$ under the operation $\sup.$
Of our four classes, $\mathcal C_1$ maps to the bottom
member of this chain, while the other three map to the top member.
Although the operation of intersection on subgroups of $S$ does
not respect the relation $\approx,$ it is not clear whether
our join-semilattice may nonetheless be a lattice.
The second author hopes to give in a forthcoming note further
results about this semilattice, and in particular, on
Question~\ref{Q.still4} above.

How much influence does the {\em isomorphism class} of a subgroup
$G\leq\Sym(\Omega)$ have on its $\approx\!\<$-equivalence class?
It does not determine that class; for
consider the abstract group $G=(\Z/p\<\Z)^\omega$ for $p$ a prime.
If for each $i\in\omega$ we let $\Sigma_i$ be a
regular $\Z/p\<\Z\!\<$-set (hence of cardinality $p)$ on which
we let $G$ act via the projection on its $i\!\<$th coordinate,
and we take for $\Omega$ a disjoint union of the $\Sigma_i,$ then
we get a representation of $G$ as a compact
subgroup of $\Sym(\Omega)$ belonging to $\mathcal{C_Q}.$

On the other hand, we may identify $G$ with the direct
product $\prod_{i\geq 0}\,(\Z/p\<\Z)^i,$ and let $\Omega$ be a disjoint
union of regular representations of the factors in this product,
getting a representation of $G$ in $\Sym(\Omega),$ also compact
in the function topology, but belonging to $\mathcal{C_P}.$
Finally, observe that if $V$ is a vector space of dimension $\aleph_0$
over the field of $p$ elements, and we also
regard $G=(\Z/p\<\Z)^\omega$ as a vector space over this field, then
$G$ and $V^\omega,$ both having the cardinality of the continuum,
are both continuum-dimensional, hence isomorphic.
Performing the same construction as before on this product expression
$G=\prod V,$ we get a representation of $G$ as a group of permutations
of a countable set $\Omega$ (with $G$ again closed in the
function topology, but no longer
compact), which Lemma~\ref{L.D_i} (with $\alpha_i$ a representative
of the $i\!\<$th orbit, and $D_i$ the product of the first $i$
orbits) shows belongs to $\mathcal C_S.$

Of course, membership of a subgroup in the class $\mathcal C_1$ is
determined by its cardinality, hence by its isomorphism class.
But to any isomorphism class $I$ of groups of continuum
cardinality, we may associate the subset of
$\{\mathcal C_S,\<\mathcal{C_P},\<\mathcal{C_Q}\}$ consisting of
those $\approx\!\<$-equivalence classes (if any) that contain
members of $I.$
Which subsets of $\{\mathcal C_S,\<\mathcal{C_P},\<\mathcal{C_Q}\}$
arise in this way (or in various related ways;
for instance, by associating to an isomorphism class $I$ the set of
{\em closed} subgroups that belonging to $I),$ we do not know.

If we take account of the topological structure of a
subgroup $G\leq S,$ this can impose restrictions on
its $\approx\!\<$-equivalence class:

\begin{lemma}\label{L.compact}
If $\Omega$ is an infinite set, then a subgroup $G\leq S=\Sym(\Omega)$
is compact in the function topology
if and only if it is closed and has finite orbits.
\end{lemma}\begin{proof}
If $G$ is closed and the members of the partition $A$ given
by the orbits of $G$ are all finite, then $G$ is a closed
subgroup of $S_{(A)}\cong\prod_{\strut\Sigma\in A}\Sym(\Sigma).$
It is not hard to see that this isomorphism is also a homeomorphism,
hence as the above product of finite discrete groups
is compact, so is $G.$

Conversely, if $G$ is compact, it is closed in $S$ by general
topology, and for each $\alpha\in\Omega$
the orbit $\alpha\<G,$ being an image of the compact group $G$
under a continuous map to the discrete space $\Omega,$ is finite.
\end{proof}

So for $|\Omega\<|=\aleph_0,$ a compact subgroup of $S$
cannot belong to $\mathcal C_S.$
Note also that if $G$ is a closed subgroup of $S$
not in $\mathcal C_S,$ then by Theorem~\ref{T.inf_orbs} there exists
a finite set $\Gamma$ such that $G_{(\Gamma)}$ has finite orbits,
so by the above lemma $G_{(\Gamma)}$ is compact.
Thus, though $G$ itself need not be compact, it will be
a countable extension of a compact subgroup that is open-closed in it,
and thus will be locally compact.

\section{Some finiteness results.}\label{S.move_fin}

This section assumes only the notation recalled in the first two
paragraphs of~\S\ref{S.Defs}, and the contents of~\S\ref{S.func_top}
(the definition of the function topology, and Lemma~\ref{L.b+f}.
At one point we will call on a result of a later section, but our
use of that result will subsequently be superseded by a more general
argument.)
We begin with a result that we will prove directly from the
definitions.

\begin{lemma}\label{L.move_fin}
Suppose $\Omega$ is a set, and $G$ a subgroup of $S=\Sym(\Omega)$ which
is discrete in the function topology on $S,$ and has the property
that each member of $G$ moves only finitely many elements of $\Omega\<.$
Then $G$ is finite.
\end{lemma}\begin{proof}
The statement that $G$ is discrete means that there is some
neighborhood of $1$ containing no other element of $G.$
Since a neighborhood basis of $1$ in $S$ is given by the
subsets $S_{(\Gamma)}$ for finite $\Gamma\subseteq\Omega\<,$
there is a finite $\Gamma$ such that $G_{(\Gamma)}=\{1\}.$

Take such a $\Gamma,$ and assuming by way of contradiction that
$G$ is infinite, let $\Gamma_0\subseteq\Gamma$ be maximal for
the property that $G_{(\Gamma_0)}$ is infinite, and let $\gamma$
be any element of $\Gamma-\Gamma_0.$
Then $G_{(\Gamma_0)}$ inherits the properties that we wish to show
lead to a contradiction; so, replacing $G$ with this subgroup,
we may assume that for some $\gamma\in\Omega,$
$G_{(\{\gamma\})},$ unlike $G,$ is finite, say of order $n.$
Then the orbit $\gamma\<G$ must be infinite, so let
$\gamma\<g_0,\ldots,\gamma\<g_n$ be $n+1$ distinct
elements of that orbit.
By hypothesis, each of $g_0,\ldots,g_n$ moves only finitely many
elements of $\Omega,$ so the infinite set $\gamma\<G$ must contain
an element $\gamma g$ not moved by any of them.
Hence $G_{(\{\gamma g\})}$ contains the $n+1$ elements $g_0,\ldots,g_n,$
contradicting the fact that, as a conjugate of $G_{(\{\gamma\})},$
it must have order $n.$
\end{proof}

If we generalize the hypothesis of this lemma by letting $\Sigma$ be
a subset of $\Omega$ and $G$ a discrete subgroup of
$\Sym(\Omega)_{\{\Sigma\}}$ each member of which moves only
finitely many elements of $\Sigma,$ it does not follow
that $G$ induces a finite subgroup of $\Sym(\Sigma).$
For example, partition an infinite set $\Omega$
into two sets $\Sigma$ and $\Omega-\Sigma$ of the same
cardinality, and let $v$ be a homomorphism
from a free group $F$ of rank $|\<\Omega\<|$ onto the group of those
permutations of $\Sigma$ that move only finitely many elements.
Take a regular representation of $F$ on $\Omega-\Sigma,$ and
consider the representation of $F$ on $\Omega=
(\Omega-\Sigma)\cup\Sigma$ given by the ``graph'' of $v,$
i.e., the set of elements of $\Sym(\Omega)_{\{\Sigma\}}$
that act on $\Omega-\Sigma$ by an element $a\in F,$
and on $\Sigma$ by $v(a).$
This subgroup is discrete because for each $\alpha\in\Omega-\Sigma$
we have $G_{(\{\alpha\})}=\{1\};$
but it does not induce a finite group of permutations on $\Sigma\<.$

On the other hand, the proof of Lemma~\ref{L.move_fin} easily
generalizes to show that if $G$ is a discrete subgroup of
$S=\Sym(\Omega),$ and $\Omega$ is the union of a family
of $G\!$-invariant subsets $\Sigma_i$ such
that every element of $G$ moves only finitely
many members of each $\Sigma_i,$ then $G$ is finite.
(Incidentally, note that throughout this section, when we refer to
families of subsets $\Sigma_i$ or $\Delta_i$ of $\Omega,$
there is no disjointness assumption.)
In a different direction, if $\Omega$ is countable we can formally
strengthen Lemma~\ref{L.move_fin} by weakening the hypothesis
``discrete'' to ``closed''; for a subgroup of $\Sym(\Omega)$ whose
members each move only finitely many elements must be
countable, and we saw in Theorem~\ref{T.countable} that a countable
closed subgroup of $\Sym(\Omega)$ is discrete.

Now suppose that for $\Omega$ countable we combine the above two
weakenings of the hypothesis of Lem\-ma~\ref{L.move_fin}, and consider
a closed subgroup $G<\Sym(\Omega)$ such that for some
expression $\Omega=\bigcup_I \Sigma_i$ of $\Omega$ as a union
of $G\!\<$-invariant subsets, each
element of $G$ moves only finitely many members of each $\Sigma_i.$
We would like to conclude that $G$ induces a finite group of
permutations of each $\Sigma_i;$ but we cannot argue as above, for now
$G$ need not be
countable, making Theorem~\ref{T.countable} inapplicable.

In an earlier version of this preprint we asked whether this conclusion
nonetheless held.
Greg Hjorth has shown us a proof, which, with his permission,
we give below.
We will use

\begin{lemma}\label{L.fin_or_inf}
Let $\Omega$ be a countable set, $G$ a closed subgroup
of $S =\Sym(\Omega),$ and $(\Delta_i)_{i\in I}$ a family of subsets
of $\Omega\<.$
Then either
\\[2pt]
{\rm(i)}~there exists a finite set $\Gamma\subseteq\Omega$ such that
$G_{(\Gamma)}\leq S_{\{\Delta_i\}}$ for all but finitely many $i\in I,$
or
\\[2pt]
{\rm(ii)}~there exists an element $g\in G$ such that
$g\notin S_{\{\Delta_i\}}$ for infinitely many $i\in I.$

Moreover, if all $\Delta_i$ are finite, then in
{\rm(i)} we can strengthen ``all but finitely many'' to ``all''.
\end{lemma}\begin{proof}
As in Lemma~\ref{L.b+f}, let $\Omega=\{\eps_0,\eps_1,\ldots\,\}.$
Assuming~(i) does not hold, we shall construct
$g_0,\,g_1,\,\ldots\in G$ which converge, by that lemma, to an
element $g$ with the property asserted in~(ii).

Suppose inductively that for some $j\geq 0$ we have chosen
$1\,{=}\,g_{-1},\,g_0,\,g_1,\,\ldots,\linebreak[0]\,g_{j-1}\in G,$
and also distinct indices $i_0,\,\ldots,\linebreak[0] i_{j-1}\in I$ and
elements $\alpha_0,\,\ldots, \alpha_{j-1}\in \Omega,$ such that for $k=
0,\ldots,j{-}1,$ $g_k$ moves $\alpha_k$ either out of or into
$\Delta_{i_k},$ and such that defining, for $0\leq k\leq j,$

\begin{xlist}\item\label{x.0-k.1,g-1}
$\Gamma_k= \{\<\eps_0,\ldots,\eps_{k-1\<}\}\,\cup\,
\{\<\eps_0\,g_{k-1}^{-1},\ldots,\eps_{k-1}\,g_{k-1}^{-1\<}\}\,\cup\,
\{\<\alpha_0,\ldots,\alpha_{k-1\<}\}\,\cup\,
\indent\{\<\alpha_0\,g_{k-1}^{-1},\ldots,\alpha_{k-1}\,g_{k-1}^{-1\<}\}$
\quad (cf.~(\ref{x.0-j.1,g-1})),
\end{xlist}
we have
\begin{xlist}\item\label{x.gkinS_*Gg}
$g_k\in G_{(\Gamma_k)}\,g_{k-1}$ for $0\leq k<j$
(cf.~(\ref{x.gjinS_*Gg})).
\end{xlist}
Since $\Gamma_j$ is finite, our assumption
that~(i) fails tells us that there are infinitely many $i$ such that
$G_{(\Gamma_j)}$ fails to preserve $\Delta_i.$
It follows that by multiplying $g_{j-1}$ on the left by a member
of $G_{(\Gamma_j)}$ if necessary, we can insure that for some index
other than $i_0,\,i_1,\,\ldots, i_{j-1},$ which we may call
$i_j,$ the resulting product $g_j$ fails to
preserve $\Delta_{i_j},$ i.e., moves an
element $\alpha_{i_j}$ into or out of $\Delta_{i_j}.$
Also,~(\ref{x.gkinS_*Gg}) shows that $g_j$ retains the
properties of the preceding elements $g_k$ of moving
$\alpha_k$ into or out of $\Delta_{i_k}$ $(0\leq k<j).$
Applying Lemma~\ref{L.b+f}, we get a limit element $g$ which
clearly preserves none of $\Delta_{i_0}, \Delta_{i_1},\ldots\,.$

To get the final assertion, observe that if all $\Delta_i$ are finite
and~(i) holds, we may take a $\Gamma$ as in~(i) and then
adjoin to it the elements of the finitely many sets $\Delta_i$
not preserved by $G_{(\Gamma)}.$
\end{proof}

Applying the above lemma (in particular
the final sentence) in the case where the $\Delta_i$ are
the singleton subsets of a set $\Sigma\subseteq\Omega,$ we get

\begin{corollary}\label{C.fin_or_inf}
Let $\Omega$ be a countable set, $G$ a closed subgroup
of $S =\Sym(\Omega),$ and $\Sigma$ a subset of $\Omega\<.$
Then either
\\[2pt]
{\rm(i)}~ there exists a finite set $\Gamma\subseteq\Omega$ such that
$G_{(\Gamma)}\leq S_{(\Sigma)},$ or
\\[2pt]
{\rm(ii)}~ there exists an element $g\in G$
which moves infinitely many members of $\Sigma\<.$\qed
\end{corollary}

We shall now get our desired result by an argument similar to
the proof of Lemma~\ref{L.move_fin}, with the above corollary replacing
our use of discreteness.

\begin{theorem}[{\rm G.~Hjorth, personal
communication}\textbf{}]\label{C.hj}
Let $\Omega$ be a countable set, $G$ a closed subgroup
of $S =\Sym(\Omega),$ and $\{\Sigma_i~|~i\in I\}$ a family
of $G\!$-invariant subsets of $\Omega$ such that each element of $G$
moves only finitely many elements of each $\Sigma_i,$
and $\bigcup_I \Sigma_i=\Omega\<.$

Then $G$ acts on each $\Sigma_i$ as a finite group of permutations;
equivalently, $G$ fixes all but a finite subset of each $\Sigma_i.$

\end{theorem}\begin{proof}
The equivalence of the two forms of the conclusion follows from the
hypothesis that each member
of $G$ moves only finitely many elements of each $\Sigma_i.$
To prove the first form of that conclusion,
suppose, on the contrary, that $G$ induces an infinite group
of permutations on $\Sigma_j$ for some $j\in I.$
Applying the preceding corollary to $\Sigma_j,$
and noting that, by hypothesis, case~(ii) of that corollary
is excluded, we get a finite $\Gamma\subseteq\Omega$ such that
$G_{(\Gamma)}\leq S_{(\Sigma_j)}.$

As in the proof of Lemma~\ref{L.move_fin},
let $\Gamma_0$ be a maximal subset of $\Gamma$ such that
$G_{(\Gamma_0)}$ induces an infinite group of permutations of
$\Sigma_j,$ and $\gamma$ any element of $\Gamma-\Gamma_0.$
Replacing $G$ by $G_{(\Gamma_0)},$ which clearly inherits the hypotheses
of the theorem, we have that $G_{(\{\gamma\})},$ unlike $G,$
induces a finite group of permutations of $\Sigma_j,$ say of order $n.$
Thus for any $g\in G,$ the group $G_{(\{\gamma g\})}$ likewise
induces a group of permutations of $\Sigma_j$ of order~$n.$

As before, $\gamma$ must have infinite orbit $\gamma\<G.$
Now applying to some $\Sigma_i$ that contains $\gamma$ the hypothesis
that each element of $G$ moves only finitely many elements of
$\Sigma_i,$ we see that each element of $G$ lies in $G_{(\{\gamma g\})}$
for all but finitely many distinct $\gamma g\in\gamma\<G.$
It follows that every finitely generated subgroup of $G$
is likewise contained in $G_{(\{\gamma g\})}$ for all but finitely many
distinct $\gamma g.$
Since $G$ induces an infinite group of permutations of $\Sigma_j,$
we can find a finitely generated subgroup of $G$ that induces a group
of $>n$ permutations of that set.
But by the last sentence of the preceding paragraph,
a group of order $>n$ can't be contained in any of the subgroups
$G_{(\{\gamma g\})},$ let alone in all but finitely many of them.
This contradiction completes the proof of the theorem.
\end{proof}

We remark that the analog of Lemma~\ref{L.move_fin} (and hence
of Theorem~\ref{C.hj}) fails for submonoids of $\Omega^\Omega\<.$
For instance, let $\Omega=\omega,$ and for
$i,j\in\omega$ define $f_i(j)=\max(i,j).$
Then $G=\{f_i\mid i\in\omega\}$ satisfies the hypotheses of
Lemma~\ref{L.move_fin} with ``submonoid
of $\omega^\omega$'' in place of ``subgroup of $\Sym(\Omega),$\!''
but does not satisfy the conclusion.

\section{Other preorderings, and
further directions for investigation.}\label{S.q.etc}

In the arguments of \S\S\ref{S.inf_orbs}-\ref{S.bdd_orbs},
when we obtained a relation $G\preq H,$
we often did this by showing that $G$ lay in the
subgroup of $S$ generated by finitely many conjugates of $H.$
This suggests

\begin{definition}\label{preqcj}
If $S$ is a group, $\kappa$ an infinite cardinal, and $G_1,$ $G_2$
subgroups of $S,$ let us write $G_1\preq_{\kappa,S}\cj G_2$
if there exists a subset $U\subseteq S$ of cardinality $<\kappa$
such that $G_1\leq \langle\,\bigcup_{f\in U}\ f^{-1}G_2 f\,\rangle.$
\end{definition}

As with $\preq_{\strt\kappa,S},$ we may omit the subscripts
$\kappa$ and $S$ from $\preq_{\kappa,S}\cj$
when their values are clear from context, and we will write
$\approx_{\kappa,S}\cj$ or
$\approx\cj$ for the induced equivalence relation.
For the remainder of this discussion, $\kappa$ will be $\aleph_0$ and
$S$ will be $\Sym(\Omega)$ for a countably infinite set $\Omega,$
and these subscripts will not be shown.

In general, $\preq\cj$ and $\approx\cj$ are finer relations
than $\preq$ and $\approx.$
Since not {\em all} the arguments in
\S\S\ref{S.inf_orbs}-\ref{S.bdd_orbs}
were based on combining conjugates of the given subgroup $G$
(in particular, some were based on conjugating a carefully constructed
element $s\in S$ {\em by} elements of $G),$ it is not obvious whether
those results can be strengthened to say that the classes of subgroups
that we proved $\approx$\!-equivalent are in fact
\mbox{$\approx\cj$\!-equivalent}.
Let us show that the answer is ``almost''.

Recall (cf.\ \cite[p.51, Theorem~6.3]{BhMMN})
that since $\Omega$ is countably infinite, the only
proper nontrivial normal subgroups of $S$ are the
group of permutations that move only finitely many points,
which we shall denote $S^{\mathrm{\<finite}},$ and the subgroup
of {\em even} permutations in $S^{\mathrm{\<finite}},$ which
we shall denote $S^{\mathrm{\<even}}.$

\begin{lemma}\label{L.when_preqcj}
Let $\Omega$ be a countably infinite set, and
$G$ a subgroup of $S=\Sym(\Omega)$ not contained
in $S^{\mathrm{\<finite}}.$
Then the unary relations $\preq\cj G$ and $\preq G$ on the set
of subgroups of $S$ coincide.
Hence if $G$ is uncountable, the unary relations
$\approx\cj G$ and $\approx G$ on that set also coincide.
\end{lemma}
\begin{proof}
In the first assertion, the nontrivial direction is
to show that $H\preq G$ implies $H\preq\cj G.$
The former condition says that $H\leq\langle\<G\,\cup U\,\rangle$ for
some finite $U\subseteq S.$
Now since $G$ is not contained in the largest proper normal subgroup
of $S,$ the normal closure of $G$ is $S.$
Hence each element of $U\subseteq S$ is a product
of finitely many conjugates of elements of $G.$
The desired conclusion follows immediately.

To see the second assertion, note that if $G$ is uncountable,
so is any $H\approx G,$ hence we also have
$H\not\subseteq S^{\mathrm{\<finite}},$ and the preceding result
can be applied to $H$ as well as to $G,$ giving $H\approx\cj G.$
\end{proof}

We can now get

\begin{proposition}\label{L.cj_eq_classes}
Let $\Omega$ be a countably infinite set, and $S=\Sym(\Omega).$
Then all $\approx$\!-equivalence classes of subgroups of $S$
other than $\mathcal C_1$ are also $\approx\cj$\!-equivalence classes.
The class $\mathcal C_1$ decomposes into the following six
$\approx\cj$\!-equivalence classes:

{\rm(i)}~The set of countable {\rm(}finite or infinite{\rm)}
subgroups not contained in $S^{\mathrm{\<finite}},$

{\rm(ii)}~The set of infinite {\rm(}necessarily countable{\rm)}
subgroups of
$S^{\mathrm{\<finite}}$ not contained in $S^{\mathrm{\<even}}.$

{\rm(iii)}~The set of infinite {\rm(}again
countable{\rm)} subgroups of $S^{\mathrm{\<even}}.$

{\rm(iv)}~The set of finite subgroups contained in
$S^{\mathrm{\<finite}}$ but not in $S^{\mathrm{\<even}}.$

{\rm(v)}~The set of finite nontrivial subgroups of
$S^{\mathrm{\<even}}.$

{\rm(vi)}~The set containing only the trivial subgroup.
\end{proposition}
\begin{proof}
The first assertion follows immediately from the preceding lemma.
In the second, it is not hard to
see that the sets (i)-(vi) partition $\mathcal C_1,$ and that
a subgroup in one of these sets cannot be $\approx\cj$\!-equivalent
to one not in that set, so it remains only to show that any two groups
in the same set in our list are $\approx\cj$\!-equivalent.

That this is true of~(i) follows from the first assertion of the
preceding lemma, and the fact that all members of $\mathcal C_1$
are $\approx$\!-equivalent.

Skipping to~(iii), if $G$ is in that class, then
Lemma~\ref{L.move_fin} shows that $G$ is non-discrete.
 From a sequence of nonidentity elements of $G$ approaching $1,$
we can extract an infinite subsequence consisting of
elements whose supports,
$\mathrm{supp}(g)=\{\alpha\in\Omega\mid g\<\alpha\neq\alpha\},$ are
pairwise disjoint.
If we take $s\in S$ whose support has singleton intersection with
each of those supports, we find that each of the corresponding
commutators $s^{-1}g^{-1}s\<g$ is a $3$\!-cycle
\cite[p.51, Exercise~6(i)]{BhMMN}.
These $3$\!-cycles lie in $\langle\<G\cup s^{-1}G\<s\<\rangle,$ and
no point belongs to the support of more than two of
them, so we can find an infinite set of
$3$\!-cycles in that group with disjoint supports.
By dropping some of these, we may assume that the complement
in $\Omega$ of the union of their supports is infinite.
Hence we may assume without loss of generality that $\Omega=\Z$ and
that we have gotten the $3$\!-cycles $(4n,\,4n\<{+}\<1,\,4n\<{+}\<2)$
for all $n\in\Z.$
Three more conjugations now give us all
$3$\!-cycles of the form $(k,\,k\<{+}\<1,\,k\<{+}\<2)$ $(k\in\Z),$
and these generate $S^{\mathrm{\<even}}.$
Hence $G\approx\cj S^{\mathrm{\<even}}.$
So all subgroups in~(iii) are $\approx\cj$\!-equivalent to that
subgroup, hence to each other, as required.

For $G$ in class~(ii), the above result shows that finitely many
conjugates of $G\cap S^{\mathrm{\<even}}$ generate
$S^{\mathrm{\<even}},$ and
since $G$ also contains an odd permutation, the corresponding conjugates
of $G$ generate $S^{\mathrm{\<finite}}.$
So all such groups are $\approx\cj$\!-equivalent to
$S^{\mathrm{\<finite}},$ and so again, to each other.
That the members of each of (iv), (v), and (vi) are mutually
$\approx\cj$\!-equivalent is easily deduced from standard results about
finite symmetric groups
\cite[\S2.4]{BhMMN}.\end{proof}

We note that the ordering on these sets induced by the
relation $\preq\cj$ on subgroups is
$$\mathrm{(i)\ \succ\cj\ (ii)\ \succ\cj\ \{(iii),\,(iv)\}
\ \succ\cj\ (v)\ \succ\cj\ (vi)},$$
with~(iii) and~(iv) incomparable.\smallskip

There is another family of preorders also implicit in the
methods we have used.
Given subgroups $G, H\leq \Sym(\Omega)$ and a cardinal $\kappa,$
let us write

\begin{xlist}\item\label{x.approx.fix}
$G\preq^{\mathrm{fix}}_\kappa H$~ if for some $\Gamma\subseteq\Omega$
with $|\<\Gamma\<|<\kappa$ we have ${G}_{(\Gamma)}\leq H,$
\end{xlist}
and let us write $G\approx^{\mathrm{fix}}_\kappa H$
for the conjunction of $G\preq^{\mathrm{fix}}_\kappa H$
and $H\preq^{\mathrm{fix}}_\kappa G.$

Lemma~\ref{L.G_*G} yields an implication between these relations and
those studied in this note:

\begin{xlist}\item\label{x.prec.emb=>prec}
$G\preq^{\mathrm{fix}}_{\strt\aleph_0}H\implies
G\preq_{\strt|\<\Omega\<|^+}H.$
\end{xlist}

The relations $\approx^{\mathrm{fix}}_\kappa$ and
$\preq^{\mathrm{fix}}_\kappa$ tend to be quite fine-grained.
For instance, given partitions $A_1$ and $A_2$ of $\Omega,$
it is not hard to see that
$S_{(A_1)}\approx^{\mathrm{fix}}_\kappa S_{(A_2)}$ if
and only if $A_1$ and $A_2$ ``disagree at $<\kappa$ elements'',
meaning that one can be obtained from the other by ``redistributing''
$<\kappa$ elements of $\Omega\<.$
\vspace{6pt}

In a different direction, one might define on abstract groups
(rather than subgroups of a fixed group) a preordering analogous
to $\preq_\kappa,$ by letting $G_1\preq^{\mathrm{emb}}_\kappa G_2$
mean that $G_1$ admits an embedding in a group $H$ which is generated
over $G_2$ by $<\kappa$ elements.
\vspace{6pt}

In our study of symmetric groups in this note,
we have considered only countable $\Omega,$ except when no additional
work or distraction was entailed by allowing greater generality.
It would be of interest to know what can be said about
$\approx_\kappa\!$-equivalence classes of closed subgroups
of $\Sym(\Omega)$ for general $\Omega$ and $\kappa\<;$ in particular,
whether there are simple criteria for a closed subgroup
$G\leq\Sym(\Omega)$ to be $\approx_{\aleph_0}\!$-equivalent
(equivalently, $\approx_{|\<\Omega\<|^+}\!$-equivalent)
to $\Sym(\Omega).$

A related topic which has been studied extensively
(e.g., \cite{SS+ST}, \cite{ST_surv})
is the {\it cofinality} of groups $\Sym(\Omega),$ defined as
the least cardinal $\kappa$ such that $\Sym(\Omega)$ can be written
as the union of a chain of $<\kappa$ proper subgroups.
If $S=\Sym(\Omega)$ is of cofinality $\geq\kappa,$
then our unary relation $\approx_{\aleph_0,\<S}S$ is equivalent to
$\approx_{\kappa,\<S}S$ (cf.\ proof of Lemma~\ref{L._0<=>_1}(ii) above);
though the converse fails under some set-theoretic assumptions.

Mesyan \cite{ZM} examines some questions similar to those considered
here for the ring of endomorphisms of the $\Omega\!$-fold direct sum
of copies of a module.

\end{document}